\documentclass[nonatbib, review]{elsarticle}
\makeatletter
\let\c@author\relax
\makeatother

\usepackage[utf8]{inputenc} 
\usepackage{hyperref}
\usepackage{lineno}
\usepackage{graphicx}
\usepackage{enumitem}
\usepackage{xspace}
\usepackage{amsmath}
\usepackage{amsfonts}
\usepackage{tikz}
\usepackage{pgfplots}
\usepackage{subcaption}
\usepackage{booktabs}
\usepackage{comment}
\usepackage{ifthen}
\usepackage[ruled]{algorithm2e}
\usepackage[backend=biber,sorting=none]{biblatex}
\usetikzlibrary{shapes, arrows.meta, calc, scopes, spy, backgrounds}
\pgfplotsset{compat=1.18}
\pgfplotsset{select coords between index/.style 2 args={
    x filter/.code={
        \ifnum\coordindex<#1\fi
        \ifnum\coordindex>#2\fi
    }
}}
\tikzset{error limit line/.style={draw=black!25, dashed, very thick}} 

\newif\ifexternalize
\newcommand{\figurefilename}[1]{\ifexternalize \tikzsetnextfilename{#1} \fi}
\ifexternalize
    \usetikzlibrary{external}
    \tikzset{external/force remake}
\fi

\definecolorset{Hsb}{wordaccent1}{}{%
a,212.46,0.6289,0.3804;%
b,212.97,0.6276,0.5686;%
c,212.73,0.5820,0.7412;%
d,212.73,0.3070,0.8431;%
e,212.73,0.1930,0.8941;%
f,212.73,0.0913,0.9451%
}
\definecolorset{Hsb}{wordaccent2}{}{%
a,0.94,0.6465,0.3882;%
b,1.25,0.6486,0.5804;%
c,1.57,0.5990,0.7529;%
d,0.87,0.3180,0.8510;%
e,1.30,0.2009,0.8980;%
f,0.0,0.0950,0.9490%
}
\definecolorset{Hsb}{wordaccent3}{}{%
a,79.66,0.5918,0.3843;%
b,79.53,0.5890,0.5725;%
c,79.59,0.5241,0.7333;%
d,80.34,0.2757,0.8392;%
e,80.0,0.1718,0.8902;%
f,81.0,0.0830,0.9451%
}
\definecolorset{Hsb}{wordaccent4}{}{%
a,267.27,0.4024,0.3216;%
b,266.94,0.4016,0.4784;%
c,267.10,0.3827,0.6353;%
d,266.84,0.1910,0.7804;%
e,268.80,0.1152,0.8510;%
f,267.69,0.0551,0.9255%
}
\definecolorset{Hsb}{wordaccent5}{}{%
a,193.52,0.6827,0.4078;%
b,193.02,0.6839,0.6078;%
c,192.68,0.6212,0.7765;%
d,192.16,0.3364,0.8627;%
e,193.20,0.2155,0.9098;%
f,192.0,0.1029,0.9529%
}
\definecolorset{Hsb}{wordaccent6}{}{%
a,27.12,0.9605,0.5961;%
b,27.10,0.9559,0.8902;%
c,27.02,0.7166,0.9686;%
d,26.92,0.4280,0.9804;%
e,27.04,0.2829,0.9843;%
f,26.67,0.1423,0.9922%
}
\pgfplotscreateplotcyclelist{mycyclelist}{%
solid,	thick,	wordaccent1a,	opacity=0.8,	mark=none\\%
densely dashed,	thick,	wordaccent1b,	opacity=0.8,	mark=none\\%
dashed,	ultra thick,	wordaccent1c,   opacity=0.8,	mark=none\\%
densely dotted,	thick,	wordaccent2a,   opacity=0.8,	mark=none\\%
dotted,	thick,	wordaccent2b,   opacity=0.8,	mark=none\\%
solid,	ultra thick,	wordaccent2c,	opacity=0.8,	mark=none\\%
}
\pgfplotscreateplotcyclelist{adaptation cyclelist}{%
solid, thick, wordaccent1a, opacity=0.8, mark=none\\%
solid, thick, wordaccent1b, opacity=0.8, mark=none\\%
solid, thick, wordaccent1c, opacity=0.8, mark=none\\%
solid, thick, wordaccent6a, opacity=0.8, mark=none\\%
solid, thick, wordaccent6b, opacity=0.8, mark=none\\%
solid, thick, wordaccent6c, opacity=0.8, mark=none\\%
}
\pgfplotscreateplotcyclelist{b spline demo}{%
solid, line width=2.2pt, black!40, mark=*, mark repeat=25, clip mode=individual\\%
solid, very thick, wordaccent2b, mark=x, mark repeat=50\\%
solid, thick, wordaccent5b, mark=o, mark repeat=50\\%
}
\pgfplotscreateplotcyclelist{complexity list}{%
solid, thick, wordaccent1a\\%
solid, thick, wordaccent1c\\%
}
\pgfplotscreateplotcyclelist{pset cyclelist}{%
draw opacity=0, mark=none\\%
error limit line, opacity=0.8, mark=none\\%
solid, thick, wordaccent1a, opacity=0.8, mark=none\\%
solid, thick, wordaccent1b, opacity=0.8, mark=none\\%
solid, thick, wordaccent1c, opacity=0.8, mark=none\\%
solid, thick, wordaccent1d, opacity=0.8, mark=none\\%
solid, thick, wordaccent1e, opacity=0.8, mark=none\\%
solid, thick, wordaccent3a, opacity=0.8, mark=none\\%
solid, thick, wordaccent3b, opacity=0.8, mark=none\\%
solid, thick, wordaccent3c, opacity=0.8, mark=none\\%
solid, thick, wordaccent4a, opacity=0.8, mark=none\\%
solid, thick, wordaccent4b, opacity=0.8, mark=none\\%
solid, thick, wordaccent4c, opacity=0.8, mark=none\\%
}
\pgfplotscreateplotcyclelist{\expandedhp cyclelist}{%
solid, thick, black, opacity=0.8, mark=none\\%
solid, thick, wordaccent1a, opacity=0.8, mark=none\\%
solid, thick, wordaccent1b, opacity=0.8, mark=none\\%
solid, thick, wordaccent1c, opacity=0.8, mark=none\\%
solid, thick, wordaccent1d, opacity=0.8, mark=none\\%
solid, thick, wordaccent1e, opacity=0.8, mark=none\\%
solid, thick, wordaccent1f, opacity=0.8, mark=none\\%
}
\pgfplotscreateplotcyclelist{\expandedts cyclelist}{%
solid, thick, black, opacity=0.8, mark=none\\%
solid, thick, wordaccent1a, opacity=0.8, mark=none\\%
solid, thick, wordaccent1b, opacity=0.8, mark=none\\%
solid, thick, wordaccent1c, opacity=0.8, mark=none\\%
solid, thick, wordaccent1d, opacity=0.8, mark=none\\%
solid, thick, wordaccent1e, opacity=0.8, mark=none\\%
solid, thick, wordaccent1f, opacity=0.8, mark=none\\%
}
\pgfplotscreateplotcyclelist{\expandedniter cyclelist}{%
solid, thick, black, opacity=0.8, mark=none\\%
solid, thick, wordaccent1a, opacity=0.8, mark=none\\%
solid, thick, wordaccent1b, opacity=0.8, mark=none\\%
solid, thick, wordaccent1c, opacity=0.8, mark=none\\%
solid, thick, wordaccent1d, opacity=0.8, mark=none\\%
}
\pgfplotscreateplotcyclelist{\expandedqs cyclelist}{%
solid, thick, black, opacity=0.8, mark=none\\%
solid, thick, wordaccent1a, opacity=0.8, mark=none\\%
solid, thick, wordaccent1b, opacity=0.8, mark=none\\%
solid, thick, wordaccent1c, opacity=0.8, mark=none\\%
solid, thick, wordaccent1d, opacity=0.8, mark=none\\%
solid, thick, wordaccent1e, opacity=0.8, mark=none\\%
solid, thick, wordaccent1f, opacity=0.8, mark=none\\%
}
\pgfplotscreateplotcyclelist{\expandedprec cyclelist}{%
solid, thick, black, opacity=0.8, mark=none\\%
solid, thick, wordaccent1a, opacity=0.8, mark=none\\%
solid, thick, wordaccent1b, opacity=0.8, mark=none\\%
solid, thick, wordaccent1c, opacity=0.8, mark=none\\%
solid, thick, wordaccent1d, opacity=0.8, mark=none\\%
solid, thick, wordaccent1e, opacity=0.8, mark=none\\%
solid, thick, wordaccent1f, opacity=0.8, mark=none\\%
}
\pgfplotscreateplotcyclelist{\expandeddecim cyclelist}{%
solid, thick, black, opacity=0.8, mark=none\\%
solid, thick, wordaccent1a, opacity=0.8, mark=none\\%
solid, thick, wordaccent1b, opacity=0.8, mark=none\\%
solid, thick, wordaccent1c, opacity=0.8, mark=none\\%
}

\title{Improving Robustness and Deployability in Distributed MPC Formation Control: a Framework for Underwater Acoustic Networks}
\author[NTNUies]{Emil~Wengle}
\author[NTNUitk]{Damiano~Varagnolo\corref{c1}}
\cortext[c1]{Corresponding author.}
\address[NTNUies]{Norwegian University of Science and Technology ({NTNU}), Department of Electronic Systems, NO-7491 Trondheim, Norway; \texttt{emil.wengle@ntnu.no}}
\address[NTNUitk]{Norwegian University of Science and Technology ({NTNU}), Department of Engineering Cybernetics, NO-7491 Trondheim, Norway; \texttt{damiano.varagnolo@ntnu.no}}
\newcommand{\niterstep}{\ensuremath{N_\text{ips}}\xspace}
\newcommand{\Ts}{\ensuremath{T_\text{s}}\xspace}
\newcommand{\Nsize}{\ensuremath{|\mathcal{V}|}\xspace}

\newcommand{\spl}{\ensuremath{\mathsf{s}}\xspace}
\newcommand{\wpl}{\ensuremath{\mathsf{w}}\xspace}

\newcommand{\zpl}{\ensuremath{\mathsf{z}}\xspace}
\newcommand{\Spl}{\ensuremath{\mathsf{S}}\xspace}

\newcommand{\dpl}{\ensuremath{\mathsf{d}}\xspace}

\newcommand{\bw}{\ensuremath{R}\xspace}

\newcommand{\psz}{\ensuremath{N_\text{pk}}\xspace}
\newcommand{\msep}{\ensuremath{\epsilon_y}\xspace}
\newcommand{\mses}{\ensuremath{\epsilon_s}\xspace}
\newcommand{\Msep}{\ensuremath{\bar{\epsilon}_y^2}\xspace}
\newcommand{\Mses}{\ensuremath{\bar{\epsilon}_s^2}\xspace}

\newcommand{\hp}{\ensuremath{h_\text{p}}\xspace}
\newcommand{\fip}{\ensuremath{w_\text{fi}}\xspace}
\newcommand{\fwidth}{\ensuremath{w_1}\xspace}

\newlength{\twocolumnwidth}
\pgfplotsset{every semilogy axis/.append style={
		width					= 0.45\twocolumnwidth,
		height					= 0.45\twocolumnwidth,
		axis on top,
		xmin					= 0,
		xmax					= 500,
		xlabel					= {time (seconds)},
		xlabel near ticks,
		xmajorgrids,
		ylabel					= {Max MSE in formation position (m\textsuperscript{2})},
		ylabel near ticks,
		ymajorgrids,
		cycle list name			= mycyclelist,
		legend pos				= outer north east, 
		legend plot pos			= left, 
		legend columns			= -1, 
		legend style			=
		{
			draw				= none,
			fill				= none,
			inner xsep			= 0.1cm,
			inner ysep			= 0.4pt,
			nodes				= {inner ysep = 0.4pt, text width = 1.5cm}, 
			cells				= {anchor = west},
		},
		legend cell align		= left, 
		grid style 				= {very thin, solid, black!20},} 
}

\begin{document}

\begin{abstract}
We consider the problem of controlling networks of underwater agents that use acoustic communication channels to exchange information, and that are tasked with maintaining a given formation while following a shared path. We thus propose approaches to robustify already existing distributed control algorithms, i.e., algorithms to transform existing distributed model predictive controllers (DMPC) unsuitable for the underwater realm (e.g., because they assume standard double precision in exchanged data, or because they are designed for synchronous, bidirectional, and reliable communication) into deployable DMPC schemes. The proposed robustification strategies accommodate both time-division and frequency-division medium access schemes, and address the inherent challenges of lossy and broadcast communication over acoustic media. More precisely, we introduce schemes for handling broadcast asynchronous communication, mitigating packet losses, and quantising exchanged data, with the goal of providing strategies that may be incorporated by other distributed formation control schemes. We characterise the data rate savings vs.\ control performance losses that may be achieved through tuning the parameters of the proposed schemes, and provide sensitivity analyses on how controller / communication performance losses are influenced by such parameters. In doing so, we showcase the implementability of the proposed strategies for practical purposes using commercially available full-duplex or half-duplex modems.
\end{abstract}


\maketitle

\section{Introduction}
\label{sec:introduction}

We consider the problem of increasing the practical usefulness of autonomous underwater vehicles (AUVs) by means of increasing the levels of distributed autonomy that they can reach. Given that using multiple vehicles is often beneficial for mission performance~\cite{Lovaas2020superUV,Sasano2016semisubmersibledev} (e.g., use multiple assets to safeguard against data loss, allow a wider search area, or reduce the deployment time of support vessels), we moreover consider distributed autonomy.

In this case communication is a means for unlocking synergies: for example, vehicles may effectively move in a planar formation~\cite{Cui10leaderfollower}, and multiple followers may maintain formation autonomously if communication is available. Simulations in~\cite{Cui10leaderfollower} show that agents indeed may keep a bounded tracking error, and perform some practically relevant distributed task -- \emph{assuming though synchronous, bidirectional, and reliable communication}.

As in large part of prior work on distributed control of AUVs, simplistic models for communication enable abstracting links as perfect information exchange along edges in an ideal static graph. A wide body of literature shows that, under ideal conditions, one may implement nominally effective distributed control schemes in underwater scenarios. But when testing in actual at-sea conditions, such algorithms may not work at all: e.g., the underwater acoustic channel is neither synchronous, bidirectional nor reliable. This suggests that the formation control algorithms that rely on simplistic communication models may fail in real-world underwater networks.

Adaptations are thus necessary, but this may introduce performance losses to the point that the envisioned distributed autonomy is practically lost. What is thus the most meaningful and generalisable way of adapting existing distributed control algorithms to make them applicable in real life situations while minimising by how much the adaptations degrade the control performance?

In this paper, we contribute towards answering this question by focusing on distributed model predictive control schemes, proposing some adaptation schemes, and assessing which performance losses such adaptations may induce. We thus present some literature in Section~\ref{ssec:literature}, before summarising our intended contributions in Section~\ref{ssec:contributions}.

\subsection{Literature review}
\label{ssec:literature}

The here cited works represent a sample of the literature on distributed control for underwater environments for which authors were successful in performing some level of field trials, or sufficiently accurate simulation environments for which the gap with actual deployment is relatively small.

Among the available works that consider the characteristics of the underwater acoustics, we consider~\cite{Chhavi20formationcontrol}, where the authors propose a gradient-descent-based delay-tolerant controller for formation control that uses a leader-follower approach. The proposed controller is tested in simulation, but also experimentally in a swimming hall, using a virtual leader. For the swimming-hall experiment, the communication from leaders to followers is expressly defined and some packets were reportedly lost, but the target path is not clearly defined, and the data -- two angles in degrees -- are sent as bit-expensive ASCII strings instead of more compact fixed-point numbers.

\cite{Zhang2023coopcontrol}, too, proposes a leader-follower cooperative control algorithm for maintaining a distance and bearing between the leader and the followers, and performs field trials for assessment purposes. Here the devices can communicate over an acoustic and an optical link. The algorithm was demonstrated to attain the target formation in a small tank; on the other hand, the authors saw the need to test the algorithm also in a less controlled environment. The optical link admits very-high-rate communication over short range and a clear channel, which is not guaranteed underwater.

\cite{elFerik13limitedBW} proposes a distributed nonlinear MPC algorithm for navigating between a sequence of points. The agents train neural networks to approximate the trajectory of the leader in interest of reducing data requirements, and make the algorithm delay-tolerant by predicting delayed packets. They demonstrate that the algorithm works under large communication delay (order of tenths of a second), which appears to account for packet duration, though it is not clear whether this is the case. Also, while the biases and weights of the neural network require only a quarter of the data of the full trajectory, the numbers themselves appear to remain floating-point numbers and invite to further compression.

\cite{Yang21realtimeDMPC} instead proposes a two-stage real-time distributed MPC algorithm with an adaptive linear discretisation approach. In the first stage, motivated by the limited communication rate available underwater, the variables and gradients are discretised to a fixed number of bits, and the discretisation interval decreases in width and centres on the new values with each iteration until an associated optimisation problem reaches a certain suboptimality. At the same time, the product of bits and iterations is bounded above by the number of bits allowed per control-loop time step. The design of the quantiser occurs offline. The second stage executes discrete-time control online, using the progressive quantiser as outlined above.

The algorithm was shown to be closed-loop stable both mathematically and by simulation, and its performance was also tested in simulation under small disturbances. Nonetheless, it remains to be seen if the algorithm would be adaptable to continuous-time systems, and the sampling time of the system appears to be very short; this suggests that packets would need to be at least as short if data are expected once every sampling period. In addition, the agents need to know everybody's initial state variable, and should a packet with quantised state be lost, the recipient has no way of recovering the lost information. 

\cite{Cai23CFC4mAUV} proposes a leader-follower approach to formation control that considers packet delays and packet loss. The controller is not based on an MPC framework, and rather a backstepping sliding-mode algorithm, which is shown to be Lyapunov convergent. Relevant to our paper is the mechanism in~\cite{Cai23CFC4mAUV} for recovering information from lost and delayed packets, that are compensated with a Gaussian predictor. The packets themselves contain state information and a timestamp, and the predictor assumes that the state variables are sampled with noise. The simulation results show that the predictor helps reducing the tracking error even under 50\% packet loss ratio. Further, the assumption of sampling state variables with noise allows applying e.g.\ quantisation noise.

\cite{DiLeo17commsaware} presents instead a communications-aware distributed MPC algorithm, in that vehicles are to maintain connectivity of the communication links, for navigating through a sequence of points while avoiding collisions. The communication is symmetric, and the agents send their planned path with a compressed format, though the compression mechanism for the data is different from the one proposed in our paper (as will be clear later, based on quantisation of opportune transformations of spline coefficients).

\cite{Hoff24commsaware} proposes a null-space-based behavioural approach to distributed formation control that is also communications-aware. The communication is defined to occur sequentially under a time-division medium access control strategy. The formation that shall be maintained is scaled to maintain a constant packet loss probability based on the experienced signal-to-noise ratio, which is modelled as a random walk in some of the case studies in the paper. The data volume to exchange is small, using six real numbers in all, and can be reduced further with knowledge of the expected range and resolution of each number.

\cite{Zhang24sparrow} proposes a formation-path-planning control algorithm that builds on the sparrow search algorithm, a swarm intelligence algorithm. The approach is leader-follower-type, in which the leader shares the trajectory with the followers and the followers coordinate their movements among themselves, and the algorithm also handles collision avoidance. The algorithm is shown to converge to formation in their case study, although communication is here assumed to be reliable (lossless), and the exact packet structure is not clear.

\cite{Gomes18attainableset} proposes an admissible-set MPC approach to control the formation of a collection of AUVs. The method here is demonstrated in simulation, and lets the vehicles avoid obstacles, although the authors assume the acoustic communication to be reliable, and the packet structure is unclear.

\cite{Hadi24reinforcement} proposes a deep reinforcement-learning approach to plan motion and control a formation. The approach is leader-follower type, where the leader is to reach a point with shortest possible safe path, and the follower is to maintain a position relative to the leader while keeping away from obstacles. Communication is considered unidirectional in one solution and bidirectional in another, and the packets appear to contain state information; information that may have compression potential. The use cases consider variable communication delay, but also reliable (lossless) links.

\subsection{Statement of contributions}
\label{ssec:contributions}

We propose a series of schemes that enable distributed MPC formation control schemes cope with the non-idealities of underwater acoustic communication (since this is the most common approach to exchange information underwater).

To exemplify how to use the proposed schemes we robustify the specific formation control algorithm in~\cite{Matous22dMPC4FPF} (henceforth referred to as ``the sensitive DMPC algorithm'' or ``the sensitive formation control scheme''). This is an alternating-direction method of multipliers (ADMM) based control scheme that assumes -- much like many of the aforementioned prior works -- synchronous, bidirectional and reliable communication of vectors of non-quantised scalars. Consequently, this scheme has unrealistically high data rate requirements for underwater acoustic standards.

We thus:
    \emph{1)} propose generalisable strategies to adapt the original scheme to the case of asynchronous, broadcast and lossy communication; 
    \emph{2)} derive an approach to compress the exchanged data to reduce the overall throughput demands;
    \emph{3)} numerically analyse the trade-offs between the adaptations and data compression schemes above, and the closed-loop performance of the derived scheme.

The objective of this paper is thus not only to propose new schemes for handling communication in formation control operations, but provide tools to better simulate which performance the agents may have in real life conditions.

\subsection{Structure of the manuscript}
\label{ssec:structure}

The remainder of the paper is structured as follows:
Section~\ref{sec:background} reviews backbone concepts and algorithms.
Section~\ref{sec:adaptations} presents and motivates the schemes made in the paper, each of which are then discussed in detail in the consequent Sections~\ref{sec:adapting-to-broadcast} (handling broadcast schemes),~\ref{sec:adapting-to-asynch} (handling asynchronous communications),~\ref{sec:handling-lost-packets} (handling packet losses), and~\ref{sec:quantizing-splines} (handling quantization of the exchanged information).
Section~\ref{sec:optimising-parameters} discusses how to optimise the parameters that define the proposed schemes,
and Section~\ref{sec:results} presents simulation results of the effects of such schemes, together with their sensitivity with respect to the parameters that define them.
Section~\ref{sec:conclusion} summarises what we believe are the main points brought by the paper.

\section{Background}
\label{sec:background}

\subsection{Essentials on Model Predictive Control (MPC)}

Referring the interested reader back to~\cite{Kouvaritakis16MPC,BorrelliFrancesco2017PCfL} for more details about MPC, this control strategy uses a (assumedly sufficiently accurate) model of the system to predict the system's states \hp steps into the future for a given control input, and uses a cost function penalizing deviation from the target state \& actuation to select, as the control inputs to be executed, the action that minimise the cost above recursively, while respecting potential constraints on both states and inputs. 

\begin{figure}
    \centering
    \begin{tikzpicture}[target pos/.style={circle, draw=black},
    good pos/.style={circle, inner sep=1.5pt, fill=black, draw=black},
    bad pos/.style={good pos, fill=red!50!black, draw=red!50!black},
    good speed/.style={-To, draw=black},
    bad speed/.style={good speed, draw=red!50!black},
    target speed/.style={good speed, help lines},
    auv/.pic={
    \draw [fill=#1, pic actions] (-0.4, -0.1) -- (-0.6, -0.2) -- (-0.6, 0.2) -- (-0.4, 0.1) -- (-0.2, 0.2)
    -- (0.4, 0.2) -- (0.6, 0.1) -- (0.6, -0.1) -- (0.4, -0.2) -- (-0.2, -0.2)
    -- cycle;
    \draw [pic actions] (-0.4, -0.1) -- (-0.4, 0.1);
    },
    auv with hand/.pic={
    \pic (-auv) at (0, 0) {auv=#1};
    \coordinate (-hand pos) at (1, 0);
    \draw [fill=#1!50!black] (-hand pos) circle [radius=3pt];
    },
    ]
    \fill [cyan!50!blue!20!white] (-1.5, -4) rectangle +(\textwidth, 0.42\textwidth);
    \pic (oops) at (0.4, 0.7) {auv with hand=red!70!yellow!50!white};
    \node (pos 0) at (1, 0) [target pos] {};
    \foreach [remember=\v as \vprev (initially 0)] \v in {1,...,3} {
        \pic (good \v) at (0, -1*\v) {auv with hand=yellow};
        \node [target pos] (pos \v) at (good \v-hand pos) {};
        \draw [help lines] (pos \vprev) -- (pos \v);
    }
    \draw [Bar-Bar] (oops-hand pos) -- (pos 0.center) node [midway, auto] {\msep};
    \pic [help lines] at (0, 0) {auv=white};
    \foreach \n in {oops, good 2} {
        \draw [good speed] (\n-hand pos) -- +(3em, 0);
    }
    \draw [target speed] (good 1-hand pos) ++ (0, -0.5em) -- +(3em, 0) node (target) [at end, inner sep=0pt] {};
    \draw [bad speed] (good 1-hand pos) ++(0, 0.5em) -- ++(4.5em, 0) node [at end, inner sep=0pt] (wrong speed) {};
    \draw [Bar-Bar] (target) -- (target -| wrong speed) node [midway, auto, swap] {\mses};
    \foreach \x/\y/\r in {0.7/-3.3/1pt,0.9/-3.45/2pt,1.15/-3.4/3pt} {
        \draw [fill=white](\x, \y) circle [radius=\r];
    }
    \draw [fill=white, decorate, decoration={bumps, amplitude=2pt}] (1.5, -3.7) rectangle (8.5, -2.3);
    \foreach [evaluate={\xactual=\x+0.1666*rand;}] \x in {1,...,5} {
        \node [target pos, inner sep=1.5pt, label=above:\x, label=below:\(\x\Ts\)] at (2+\xactual, -3) {};
        \draw (\x+2, -3) circle [radius=0.4pt];
    }
    \node [target pos, inner sep=2pt, label=above:here, label=below:\phantom{I}now\phantom{I}] at (2, -3) {};
    \draw (2, -3) circle [radius=0.4pt];
    \node [target pos, inner sep=2pt, label=above:\smash{\hp}, label=below:\(\hp\Ts\)] at (8, -3) {};
    \draw (8, -3) circle [radius=0.4pt];
    \node at (5.5, -2) {the AUV predicts\dots};
\end{tikzpicture}  
    \caption{Illustrating some key concepts in this paper.
    }
    \label{fig:formation-error}
\end{figure}
Fig.~\ref{fig:formation-error} illustrates what is meant by the deviation from the target state. Four autonomous underwater vehicles are travelling forward in a wide formation, using their thrusters and fins (their actuation inputs) to keep their \emph{hand position} centred in the large circles in front of them. The topmost vehicle has come off course, so its hand position has a non-zero position error \msep. The second vehicle is travelling faster in the forward direction than the others, which results in a non-zero speed error. The bottommost vehicle is thinking of where it predicts to be in the near future, up to \hp steps ahead. Each step represents a sample \Ts ahead of the previous one, starting from the present. This forecast lets it adjust its planned actuation so it can stay closer to its objective of maintaining constant speed, which would lead to a prediction aligned with the tiny circles.

Distributed MPC extends MPC to networks of systems by defining an opportune ``network-oriented'' control problem, and by letting the agents solve it in a collaborative manner~\cite{Camponogara02DMPC}. Typically the collaborating agents do not require global information (e.g., need not know the state of everybody else) but rather have access to the own state and partial information about their neighbours, and ladder on opportune distributed optimization methods.

Consider then that our quest is for general strategies for dealing with communication non-idealities; at the same time, we shall demonstrate them on some specific formation control algorithm.
Among the several approaches to implement distributed MPC in the literature~\cite{Gruene2017nmpc,Papaioannou23targetjamming}, we thus focus on the specific algorithm~\cite{Matous22dMPC4FPF}, that is based on the rather general and widely used Alternating Direction Method of Multipliers (ADMM) numerical optimisation algorithm, and that focuses on the generic task of formation path following (\emph{1)} letting a collection of underwater vehicles maintain formation, and \emph{2)} letting the barycentre of the formation follow a path at a constant speed). By focusing on this scheme we thus consider a situation that several other authors have considered before~(\cite{saska2016predictive,vanParys17dMPC4MVS}). We nonetheless briefly describe this specific controller below, so to make the manuscript sufficiently self-contained.

\subsection{Essentials on~\cite{Matous22dMPC4FPF}}

The objective function proposed in~\cite{Matous22dMPC4FPF} is
\begin{equation}
    \label{eq:matous-dmpc}
    J_{0,i}
    =
    \int\limits_t^{t+T}
    \left(
    {\widetilde{p}_i^\mathrm{T}(\tau)
    Q_p
    \widetilde{p}_i(\tau)}
    +
    Q_s
    {\left(\dot{s}_i(\tau)-U_d\right)}^2
    \right)
    \mathrm{d}\tau,
\end{equation}
where \(\widetilde{p}_i\) is the deviation of vehicle \(i\)'s position \(p_i\) from its desired position \(\widehat{p}_i\) in the formation,
\(Q_p\) is its associated positive semidefinite penalty matrix,
\(\dot{s}_i\) is the speed of the vehicle along the predefined path,
\(U_d\) is the desired speed along the path, and
\(Q_s\) is a non-negative penalty factor for speed error.
Since \(s_i(0)=0\), \(s_i(t)\) is the distance vehicle \(i\) has travelled along the target path by the time \(t\).

When optimised in \(p_i(t)\) and \(s_i(t)\) over the integration interval,~\eqref{eq:matous-dmpc} is a variational problem. To reduce it to a finite dimensional problem, \(p_i(t)\) and \(s_i(t)\) are restricted to be spline functions, so that the optimisation vector becomes the spline coefficients.
This spline approximation technique imposes restrictive constraints on the control systems, however:
the models must be differentially flat, which means that their state \(x\) (\(p\) being part of the state here) and actuation \(u\) must be obtainable from their output \(y\), its derivatives and its antiderivatives.

In this work, the vehicles are assumed to be six-degrees-of-freedom (6DOF) marine vehicles (a choice reflecting the one made in~\cite{Fossen11marinecraft} and evaluated in~\cite[c.\ 9]{Matous23thesis}).
The vehicles are further assumed to be underactuated, and the controller to operate on a hand position, a point that lies at a positive distance in the vehicle's forward direction from its barycentre (i.e. the point in front of the AUVs in Fig.~\ref{fig:formation-error}; a technical choice that allows an opportune output-feedback linearisation of the model, transforming it to a double integrator).

Combining this opportune feedback linearization approach with the differentially flat structure of the model, \eqref{eq:matous-dmpc} reformulates the MPC optimisation problem so that its search space is that of the spline coefficients that represent the own positions \(p\), and the relative positions of the other agents represented by \(s\).

We note that this is a rather classical trajectory-generation approach, and several other algorithms rely on flatness (and thus on polynomial/spline trajectories). For vehicles that have slow dynamics, smooth trajectories, significant communication constraints, and moderate required accelerations. In this scenario flatness-based spline parameterisations work well and do not limit performance (while for highly dynamic maneuvers like high-speed marine interaction with waves this may not be the case). For example,~\cite{vanParys17dMPC4MVS, Dinh24uavswarm, Mulagaleti18flexpayload} too base their formation controllers on splines / polynomial coefficients.

\subsection{Essentials on the Alternating Direction Method of Multipliers (ADMM)}

ADMM solves optimisation problems in \(x, z\) on the form \(J_0(x, z)=f(x)+g(z)\), subject to the affine equality constraint \(Ax+Bz=c\). In each iteration, ADMM optimises \(J_0\) first in \(x\), then in \(z\), and lastly updates the dual multipliers associated with the equality constraints. ADMM typically finds a solution with ``modest accuracy (\dots) within a few tens of iterations''~\cite{Boyd11ADMM}, and can also be made distributed, as is the case in numerous applications~\cite{Boyd11ADMM,vanParys17dMPC4MVS,Bastianello21asyncADMM}.

In its basic form, ADMM requires some level of coordination among the agents in terms of which agent communicates first, etc. This may be obtained in different ways, ``Time Division Multiple Access'' and ``Frequency Division Multiple Access'' being likely the most common ones. The next two subsections introduce them below, again to keep the paper self-contained.

\subsection{Essentials on Time Division Multiple Access (TDMA)}
\label{sub:tdma}

TDMA is a Medium Access Control (MAC) strategy that divides the digital channel in non-overlapping and cyclically repeating time slots. Users are then assigned slots where they are allowed to transmit data~\cite{Goldsmith05wireless} (see Fig.~\ref{fig:tdma-diagram} for a sketch). Each block shows which user is allowed to transmit where.

\begin{figure}[!htbp]
    \centering
    \hfill\begin{subfigure}[t]{0.5\twocolumnwidth}
        \begin{tikzpicture}[
time slot/.style={rectangle, draw, very thick, minimum width=3.15cm, minimum height=0.4cm, text width=2.5cm, align=center, anchor=south west, font=\small,
 inner sep=1pt},
]
    \draw [To-To] (4, 0) -- ++ (-4, 0) node [below, near start] {frequency} -- ++ (0, 4) node [right] {time};
    \foreach [evaluate={\usr=int(mod(\x,4)); \ypos=\x*0.5;}] \x in {0,...,5} {
        \node [time slot] at (0, \ypos) {User~\usr};
    }
\end{tikzpicture}
        \caption{TDMA.}
        \label{fig:tdma-diagram}
    \end{subfigure}\hfill%
    \begin{subfigure}[t]{0.5\twocolumnwidth}
        \begin{tikzpicture}[
time slot/.style={rectangle, draw, very thick, minimum width=3.15cm, minimum height=0.6875cm, text width=2.5cm, align=center, anchor=south west, font=\small,
 inner sep=1pt},
time slot/.default={},
]
    \draw [To-To] (4, 0) -- ++ (-4, 0) node [below, near start] {frequency} -- ++ (0, 4) node [right] {time};
    {[rotate=-90, transform shape]
    \foreach [evaluate={\usr=int(mod(\x,4)); \ypos=0.1+\x*0.7875;}] \x in {0,...,3} {
        \node [time slot] at (-3.15, \ypos) {User~\usr};
    }
    }
\end{tikzpicture}
        \caption{FDMA.}
        \label{fig:fdma-diagram}
    \end{subfigure}\hfill%
    \caption{Representation of how the channel is divided in the considered MAC protocols. Note that guard periods and/or guard bands may be added to mitigate the risk of packet collisions and interference, trading off spectral efficiency against robustness.}
    \label{fig:mac-diagrams}
\end{figure}

TDMA requires a common perception of the system time across the network to be effective. Indeed, clock offset and skew -- if present -- must be compensated, lest the performance of the MAC scheme eventually degrades when the resulting drift becomes too severe. Distributed clock synchronisation algorithms, examples of which are given in~\cite{Bolognani09PIclocksync,Liu14dasync,Bolognani16randomclocksync}, can be executed in parallel with the control ones considered in this paper to attain a sufficiently accurate synchronisation. We note that the communication overhead coming from implementing these synchronisation algorithms is quite limited (for example, in~\cite{Bolognani09PIclocksync}, where the agents send only one number, such requirements were in the worst case 64~bits for guaranteeing the maximum deviation among the various clocks of a 10-agent network to be smaller than 10~ms in 200~iterations). For this reason, we neglect this requirement in this paper and assume the clocks to be sufficiently synchronised to enable implementing TDMA effectively.

We also note that underwater wireless networks often suffer from significant propagation delay not existing with radio communication (the nominal sound speed in water is 1500~m/s), so network protocols that solve terrestrial-network problems are often of little use underwater~\cite{editorialJOE19uac}. 
Variants of TDMA exist for networks with non-negligible latency, however. Rather than considering the propagation delay as a burden, the work in~\cite{Chitre12delaytolerant} shows that delay-tolerant MAC protocols can exploit it to raise the network-wide throughput. One example of such a protocol is super-TDMA~\cite{Lmai2017superTDMA}, which allows multiple agents to transmit packets at once. All the use cases in the paper consider the links to be short, so the propagation delay is on the order of single tens of milliseconds, which we believe is short enough to make TDMA practically usable nonetheless. 

\subsection{Essentials on Frequency Division Multiple Access}
\label{sub:fdma}

Frequency Division Multiple Access (FDMA) is a similar MAC strategy to TDMA. The difference is that the channel is divided in non-overlapping frequency bands rather than in cyclically repeating time slots. FDMA can be implemented with guard bands to protect users from interference from Doppler spread in the channel. The frequency bands are assigned to users~\cite{Goldsmith05wireless}, so FDMA enables concurrent communication at the expense of decreased per-user bandwidth (see  Fig.~\ref{fig:fdma-diagram} for a sketch).

Using FDMA in a broadcast-based network effectively requires full-duplex capabilities, which means the users must be capable of transmitting and receiving simultaneously (otherwise agents that are broadcasting would miss all incoming packets, or be unable to transmit at all because always busy receiving packets).

FDMA by itself does not require clock synchronisation to be effective, so the agents may update asynchronously in an FDMA setup. However, full-duplex FDMA supports synchronous communication. As outlined in Section~\ref{sub:tdma}, synchronous communication depends on a distributed clock synchronisation algorithm whose payload can be piggybacked on the DMPC packets. We therefore assume that the clocks are synced.

\section{On the need of adapting DMPCs to real world underwater acoustic communication}
\label{sec:adaptations}

Several DMPCs assume agents to communicate with a synchronous, bidirectional, and reliable channel. Underwater acoustic communication is typically asynchronous, unidirectional (broadcast), and unreliable -- the opposite of the current assumptions. Packets are broadcast, propagation delay is non-negligible, and the channel is highly diverse -- see for example~\cite{vanWalree13propagation} for a comprehensive study on channel diversity. Without adopting a medium access control (MAC) technique and/or enabling full-duplex operation (simultaneous transmission and reception), packets will collide and be lost, or at least interfere with each other. 

Such DMPCs need a few adaptations if we wish to use them for practical acoustics-based formation control of underwater vehicles. Hence, the following sections list and discuss such adaptations.

\section{Adapting to broadcast communication}
\label{sec:adapting-to-broadcast}

The scheme below applies to algorithms that suppose each agent to perform \niterstep local iterations of updating the control inputs per prediction interval \Ts and then exchanging such results. In this case the original algorithm would assume that \(\Nsize(\Nsize-1)\) packets are sent simultaneously over multiple communication links, where \(\Nsize\) is the number of vehicles in the network.

Very often controllers are implemented so that every agent exchanges one packet with each neighbour at the end of each ADMM iteration. As argued before, underwater acoustic communication is inherently broadcast rather than unicast. Even if one node is the intended recipient of a packet, every node within reception range could potentially receive the same packet. Given this inherently broadcast nature of the medium, the first need is to modify the structure of the algorithms to support broadcast communication.

The easiest method is to combine data that would otherwise be sent in one packet per neighbour in a single packet (i.e., the agent ``sending out'' its information sends everything it needs to communicate into one packet to be broadcast, and structured in a way that the receiving neighbours can find the information intended for them within this packet by opportunely looking at a specific part of the message).

More explicitly, if we let \(\wpl_{BA}\) be the generic vector of scalar values sent from $A$ to $B$ in the original DMPC algorithm, then concatenating all the \(\wpl_{BA}\)'s for the various $B$'s into the payload of $A$, and appending a numeric identifier to such a payload, then agent~\(A\) can send all necessary information to its peers in one packet per ADMM iteration. Note that this labelling may be tacit if the topology of the network is guaranteed to be static.

We note that other broadcast-based ADMM strategies may be implemented. For example, D-ADMM~\cite{Mota13dADMMseparable,Mota15dADMMlocaldomain} lets each agent in the network keep a local copy of the optimisation variable \(x_i\) and broadcast it to its neighbours. The agents send their updates in a sequence determined by a graph colouring scheme to avoid packet collisions, and may allow multiple agents to transmit at once if they are sufficiently far apart. In cases where networks that have a small graph diameter, running a graph colouring algorithm would only serve to produce unnecessary overhead.

\section{Adapting to asynchronous communication}
\label{sec:adapting-to-asynch}

When an algorithm assumes that all communication occurs simultaneously and synchronously, it needs to be modified so to avoid packet collisions. To aid this we propose two alternative strategies:
\begin{enumerate}
    \item communicate asynchronously, by assigning each agent its own frequency band in an FDMA scheme, to be used to send the own information to the peers;
    \item communicate sequentially, by implementing a TDMA scheme, thus defining implicitly who communicates when.
\end{enumerate}
The two strategies are dual with respect to each other: either separate the agents in time or separate them in frequency. The former requires keeping clocks synchronised in a distributed way (that may be implemented, e.g., via~\cite{Bolognani16randomclocksync,Liu14dasync} or similar approaches, at the cost of some small communication overhead that is ignored in this paper).

As for the FDMA strategy, the agents may parse information sent from different agents by opportunely processing the recordings in separate frequency bands, using the bands to sort the information by sender agent. This strategy implies that packets may potentially be sent simultaneously, and this alleviates the need for sequential communication. At the same time, using FDMA effectively requires full-duplex capabilities, increases the usage of the acoustic spectrum, and reduces the available bandwidth per transmitting agent by a factor \Nsize. The reduction would in practice be somewhat bigger due to the need for guard bands to defend against Doppler spread. In return, the agents would know whose frequency band is whose, so there would be no need for a numeric identifier in the payload.

As for the TDMA strategy, the agents run iterations sequentially in a round-robin fashion, such that one agent takes an update step, broadcasts its packet, then the next agent takes these actions, to avoid packet collisions. This situation requires more strict assumptions on the clock synchronisation of the agents, and also implies a factor \Nsize less available time per iteration than the FDMA approach, though it does not constrain the available frequency band. In TDMA settings, the transmitting agent needs to include its identifier in the packet it sends, and assuming that the sender knows who is in the network, the recipients can determine which part of the payload is theirs. Under an effective clock-synchronisation scheme, the selection of which is beyond the scope of this paper, scheduling by TDMA can be used to implement the sequential iterations and broadcasts. A timestamp can also be sent piggyback on the payload in a practical implementation.

Worth noting is that the \Nsize is, for both FDMA and TDMA, an upper bound on the factor for which we must divide the time or frequency budget per agent, neglecting the effect of guard regions. Agents that are at least three communication hops apart can potentially use the same time slot or frequency band without causing interference. Determining the actual factor \(\Nsize'\) is a 2-hop graph colouring problem, for which there exist fast solution algorithms, e.g.,~\cite{Schneider2009graphcolouring,Ahmed2010distributedgraphcolouring}. The data rate requirements we give are therefore upper bounds, tight for networks for which there exists no shortest path of three hops or longer.

\section{Adapting to lossy channels}
\label{sec:handling-lost-packets}

Assuming reliable communication in underwater realms is risky, since in this case packets may be lost (and algorithms may fail to ensure convergence). The mechanism we propose to handle unreliable communication is to suppose that agents can infer whether a specific packet has been lost or not. Assuming that the network employs a simple TDMA scheduler (and similar considerations may be done in the frequency domain) actually enables this: if agent $B$ knows that agent $A$ should have sent it a packet in the near past, but $B$ has not received this expected packet from $A$ in its expected time slot, then $B$ may infer that the packet from $A$ has been lost.

\begin{figure}[!htbp]
    \centering
    \begin{tikzpicture}[x=1.5em, y=1.5em]
\def\lawnmowerwidth{10}
\def\lawnmowerradius{1}
    \draw [-Stealth, thick] (0, 0) -- ++ (\lawnmowerwidth, 0) arc (-90:90:\lawnmowerradius) -- ++ (-\lawnmowerwidth, 0) arc (270:90:\lawnmowerradius) {[dashed] -- ++ (\lawnmowerwidth, 0)};
    \draw [help lines, decorate, decoration=brace, thick] (0, 0.1) -- ++ (\lawnmowerwidth, 0) node [midway, above] {\(l\)};
    \draw [help lines] (\lawnmowerwidth, \lawnmowerradius) -- ++ (-30:\lawnmowerradius) node [midway, above] {\(r\)};
\end{tikzpicture}
    \caption{A schematic representation of a lawnmower path, along with its defining parameters. Agents are supposed to follow this path using a constant speed.}
    \label{fig:lawnmower}
\end{figure}
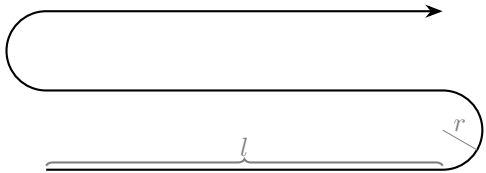

To handle this event, one may develop ad-hoc data imputation algorithms to reconstruct the lost information. Note that in this step there is the need for ad-hoc considerations, since different information contents may benefit from different data imputation approaches. As an example (that, however, may not be suitable for network control of vehicles with fast and unstable dynamics), consider algorithms where agents implicitly know their target paths / target trajectories (as in~\cite{Matous22dMPC4FPF} and in several other related papers). For example, consider a lawnmower path as in Fig.~\ref{fig:lawnmower}, typically used in seabed mapping or in localisation scenarios, where agents may be seeking a specific object for recovery or destruction purposes. Here, a given number of vessels shall ideally keep a line formation perpendicular to the straight sections of the lawnmower path (thus not turning with the path), and distributed such that the barycentre of the vehicles is on this path.

In such a scenario, each agent knows the nominal target path every agent should follow in space. The exchanged packets may thus encode forecasts of how far the agents believe the vehicles will have travelled in the near future, and this information is used in the DMPC to update the agents' control inputs to keep the formation as close to the desired one as possible.

Formally, the information encoded by the generic agent $A$ in one of its packets may be represented via the equation
\begin{equation}
    \label{eq:packets}
    \wpl_{BA}
    = 
    2\rho\spl_A
    -
    \zpl_{BA},
\end{equation}
where
\begin{itemize}
    \item \(\rho\) is a real scalar value, representing a parameter proper of the ADMM algorithm and not of interest in this paper,
    \item \(\spl_A\) is a vector of real values, representing how far vehicle~\(A\) plans to have travelled along the path it is tasked to follow in a the next time window, and
    \item \(\zpl_{BA}\) is a vector of real values, representing what agent \(A\) believes how far vehicle~\(B\) plans to have travelled along the path it is tasked to follow in a the next time window (in other words, \(\zpl_{BA}\) is what \(A\) believes \(B\) has as its own \(\spl_B\)).
\end{itemize}

These forecasts of how far an agent is along its target path are encoded as uniform B-splines, i.e., as linear combinations of piecewise smooth third-degree polynomials over a span of a finite number of equally spaced break points in time (intuitively, the unique points in time when third derivatives are allowed to change; two consecutive break points define an \emph{interval} of such a spline). Each polynomial in the linear combinations above spans up to four contiguous intervals, is zero outside these intervals, and is opportunely normalised. 

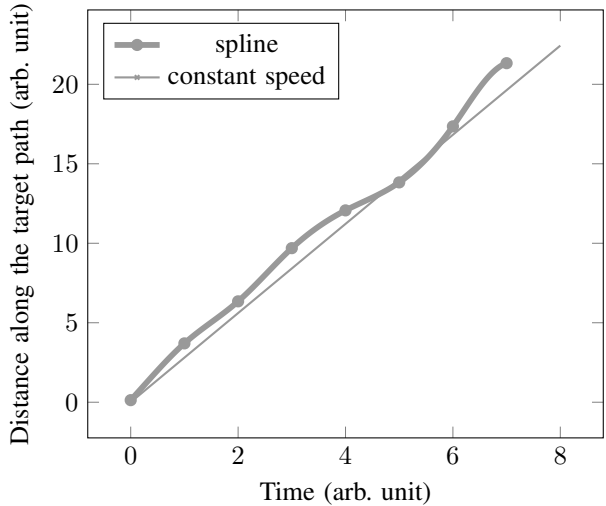
\begin{figure}[!htbp]
    \centering
    \figurefilename{tikz-b-spline-first}

\tikzset{every mark/.append style={scale=0.6},}
\begin{tikzpicture}
        \begin{axis}
        [
            xlabel                  = {Time (arb. unit)},
            ylabel                  = {Distance along the target path (arb. unit)},
            cycle list name         = b spline demo,
            mark size               = 2pt,
            legend pos              = north west,
        ]
               	\addplot+ [select coords between index={0}{350}, mark repeat=50] table [ x = t, y = target ]%
    			{Plotdata_b-spline_extrapolate-for-real.txt};%
    			\addlegendentry{spline}
            \addplot+ [black!40, thick] coordinates {(0,0) (8,22.434668)};
            \addlegendentry{constant speed};
        \end{axis}
\end{tikzpicture}
    \caption{A graphical visualisation of the information encoded in a packet. The dots along the blue curve correspond to the break points, and thus indirectly give information on the tangential speed that the vehicle is supposed to follow in time. For reference, a straight line would indicate a constant speed along the planned trajectory.}
    \label{fig:prev-packet}
\end{figure}

Thus, these splines represent how far a vehicle is supposed to be along a known trajectory in space, and thus positions in space at a certain time. Fig.~\ref{fig:prev-packet} depicts how one may represent the evolution along trajectories with splines.

Assume that agent~\(B\) does not receive the information~\eqref{eq:packets} from agent~\(A\). This means that \(B\) does not have the most recent update of where \(A\) forecasts itself to be at a given time, and of where \(A\) believes \(B\) will be too. A statistically meaningful way for \(B\) to impute the beliefs of \(A\) is to opportunely extend the older beliefs of \(A\), i.e., extrapolate the splines encoded in the latest packets it received from \(A\) beyond their original time horizon. In the next two subsections, we present two candidate methods of extrapolation.

\subsection{Extrapolating jerk}
\label{ssec:naive-extrapolation}

One may extrapolate a position-along-a-trajectory by assuming that the vehicle will keep the same jerk (the first time derivative of acceleration) after the end of the original time horizon.

This may be done mathematically by adding one break point (and thus an interval) at the endpoint of the original spline and keep the same jerk for the length of that added interval. A positive side of letting the B-spline continue this way is that it does not affect the shape of the spline in the original intervals, and this means not changing the original forecasted distances along the path in the original time domain.

However, if the acceleration or the jerk is nonzero at the endpoint of the original trajectory (as is the case whenever the velocity of the agent is forecasted to not be constant over the last interval), then the extrapolated interval may cause a forecasted distance along the path that substantially diverges according to the last acceleration and jerk. This effect can be seen graphically in Fig.~\ref{fig:extrapol}.

\subsection{Extrapolating velocity}
\label{ssec:fitting-extrapolation}

Alternatively, one may extrapolate a distance along the path assuming that the vehicle will continue after the end of the original time horizon at the same velocity, by artificially setting it to have zero acceleration and zero jerk.

Extending in this way means in general not being able to encode this distance along the path information as a uniform B-spline, since forcing the acceleration to zero may introduce a discontinuity in that signal. To be compliant with the original purpose this information shall serve, there is the need to perform an opportune re-sampling of the original spline.

To do so, we can let the to-be-extended B-spline \(\Spl(t)\) be defined by the temporal break points \([0,1,\dotsc,k]\), such that \Spl spans the interval \(t\in[0, k]\). One can then evaluate the original \Spl at regular temporal instants \([0,1/m,\dotsc, k]\) first, and then append a sequence of future (fitting) points \([k+1/m,k+2/m,\dotsc,k+1]\) to which one can assign values such that the derivative with respect to \(t\) is fixed as the left derivative in \(t=k\). Finally, to construct the extended B-spline over the interval \([0, k+1]\), we may find the coefficients of a B-spline with interval \([0, k+1]\) and break points \([0, 1,\dotsc, k+1]\) using the sample points as fitting data and use some fitting method for the purpose (in this paper, a simple least squares approach and letting \(m=2\) -- that is, sample twice per interval -- was enough to obtain practically useful results).

The extrapolating velocity mechanism leads to a trajectory that matches the first left derivative at the end of the original spline, but the zero second derivative introduces a discontinuity; this leads to a recomputing step that introduces some fitting errors almost everywhere in the new spline. In other words, the original forecasted distance along the path in the original time domain is perturbed.

\subsection{Comparing the two extrapolation methods}

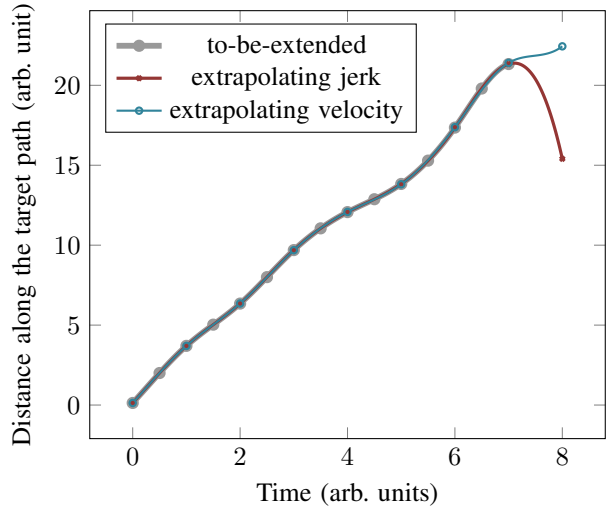
\begin{figure}[!htbp]
    \centering
    \figurefilename{tikz-b-spline-demo}

\tikzset{every mark/.append style={scale=0.6},}
\begin{tikzpicture}
        \begin{axis}
        [
            xlabel                  = {Time (arb. units)},
            ylabel                  = {Distance along the target path (arb. unit)},
            cycle list name         = b spline demo,
            mark size               = 2pt,
            legend pos              = north west
        ]
            \foreach \mypath/\mypathlabel/\lastdatapoint in {target/to-be-extended/350,naive/extrapolating jerk/400,linear/extrapolating velocity/400}
                {
                \edef\addseries{\noexpand%
               	\addplot+ [select coords between index={0}{\lastdatapoint}] table [ x = t, y = \mypath ]%
    			{Plotdata_b-spline_extrapolate-for-real.txt};%
    			\noexpand\addlegendentry{\mypathlabel};
                }
                \addseries
    		}
	    \end{axis}
\end{tikzpicture}
    \caption{A graphical comparison of the two extrapolation mechanisms described in Sections~\ref{ssec:naive-extrapolation} and~\ref{ssec:fitting-extrapolation}. The markers in the to-be-extended spline are located at the points at which the original spline is sampled for the extrapolating velocity approach. Every other marker is instead a break point for the associated uniform B-spline.}
    \label{fig:extrapol}
\end{figure}

As Fig.~\ref{fig:extrapol} shows, the jerk-extrapolation approach may lead to, in a sense, aggressive extensions that may significantly over- or under-estimate where the agents will be at the end of the extended time horizon. This may then induce the distributed MPC scheme to increase the actuation signals significantly (i.e., make the agents more reactive than necessary). 

\begin{figure}[!htbp]
    \centering
    \tikzset{every mark/.append style={scale=0.6},}
\begin{tikzpicture}
        \begin{axis}
        [
            xlabel                  = {Time (arb. units)},
            ylabel                  = {Distance along the target path (arb. unit)},
            cycle list name         = b spline demo,
            mark size               = 2pt,
            xmin=6.0, xmax=7.3,
            legend pos              = south east,
        ]
            \foreach \mypath/\mypathlabel/\lastdatapoint in {target/to-be-extended/350,naive/extrapolating jerk/400,linear/extrapolating velocity/400}
                {
                \edef\addseries{\noexpand%
               	\addplot+ [select coords between index={0}{\lastdatapoint}] table [ x = t, y = \mypath ]%
    			{Plotdata_b-spline_extrapolate-for-real.txt};%
    			\noexpand\addlegendentry{\mypathlabel};
                }
                \addseries
    		}
	    \end{axis}
\end{tikzpicture}
    \caption{Comparing the two extrapolation methods against each other around the end of the old time horizon. Extrapolating jerk does not affect the existing forecast, while extrapolating velocity introduces some perturbation towards the end of it.}
    \label{fig:extrapol-zoom}
\end{figure}
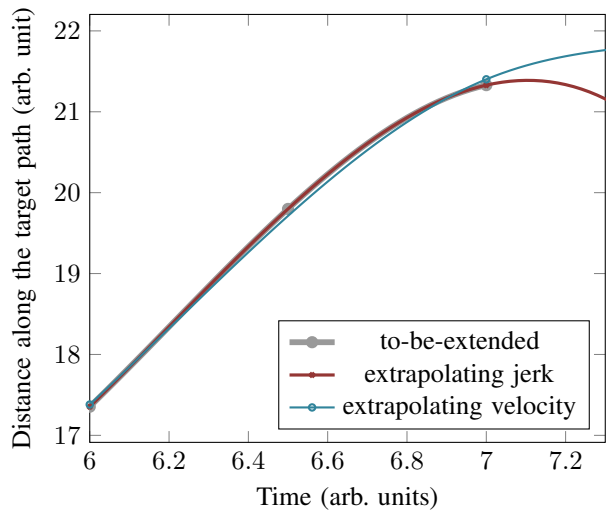

The velocity-extrapolation scheme may instead introduce some unwanted ripples near the end of the old time horizon, as Fig.~\ref{fig:extrapol-zoom} highlights. These ripples are the consequence of approximating the discontinuity in second derivative at the point where the extended trajectory begins. This may too be accounted for by the distributed MPC scheme by increasing the actuation signals (i.e., again make the agents more reactive than necessary).

The two schemes thus present this trade-off:
\begin{description}
    \item[extrapolating jerk] does not perturb the final part of the original forecast, but may over-extend it in the wrong way; 
    \item[extrapolating velocity] statistically does not present over-extensions, but may distort the forecast towards the end of the original time horizon.
\end{description}

Our simulations-based analyses strongly suggest that, at least for the two use cases we consider in Section~\ref{sec:results} the best approach in practical scenarios is to extrapolate velocities. So, even if managing packet-losses by means of the two data-imputation techniques proposed here means translating the loss of a packet into an opportune distortion of the final temporal parts of the forecasts exchanged by the agents, extrapolating velocity seems to always guarantee a better closed loop performance. Hence, Section~\ref{sec:results} will only present results relying on the velocity-based approach.

\subsection{Additional factors aiding the practical implementation of the scheme}

Any extrapolation method shall preserve the temporal horizon of the predicted trajectories. This may be done by effectively discarding the oldest \(\tau\) seconds of the prediction one would get by means of the extrapolation mechanisms above. This is practically done by transforming the coefficients of the extrapolated spline with an opportune basis transformation~\cite{vanParys17dMPC4MVS} into a spline with domain \([\tau, \hp\Ts+\tau]\) and the same number of intervals as the original spline, where \(\tau\) is the desired prediction time increment.

Note that, often, formation control DMPC algorithms use splines extrapolation on the results from one iteration of ADMM to obtain a better starting point for the next ADMM iteration. The motivation behind such an extrapolation is that the next ADMM iteration is executed a non-negligible amount of time in the future. The extrapolation hence ``warm-starts''~\cite{vanParys17dMPC4MVS} the optimisation problem~\eqref{eq:matous-dmpc}.
Note also that a similar temporal shifting mechanism may be implemented also when successfully receiving the packets, but with a latency that cannot be neglected. Indeed, with non-negligible delay, packets will not contain up-to-date samples of \(\wpl_{BA}(t)\) (i.e., forecasted trajectories of the agents), but samples of \(\wpl_{BA}(\tau_A)\) at a time instant \(\tau_A<t\). This typically implies that each agent believes that they are ahead of the rest of the network, to slow down as a consequence, and possibly cause the formation to halt. In this case, a trajectory-extrapolation step based on the extrapolation used to compensate for losing a communication packet may compensate such a lag \(t-\tau_A\).

\section{Quantising to to-be-exchanged splines coefficients}
\label{sec:quantizing-splines}

Formation control DMPC algorithms for the underwater realms use very often double-precision floating-point numbers in all communication, meaning that if we were to represent 1~km this way, the double-precision representation would have significant figures down to the picometre. Relaxing the precision requirements enables compressing the data, which is instrumental to save data rate.

Commercially available underwater acoustic modems can output between tens of bits per second and tens of kilobits per second~\cite{Kilfoyle00,Zia21modems}, making data compression paramount to enable implementation on more types of modems. Even under significant data compression, packets could last for hundreds of milliseconds, which limits the packet transmission rate. The limited packet rate is further restricted by the non-negligible propagation delay, because the nominal speed of sound in water is near 1500~m/s.

The data compression strategy we propose is that one that numerically seems to have the best performance among the different options we investigated, and it ladders on the following intuition: when a formation control DMPC algorithm sends vectors of numbers that encode opportune splines representing predicted positions in the near future, we may compress the information content by identifying the typical behaviour encoded by such splines, and then exploit it to find an opportune encoding (and thus send opportunely approximated deviations).

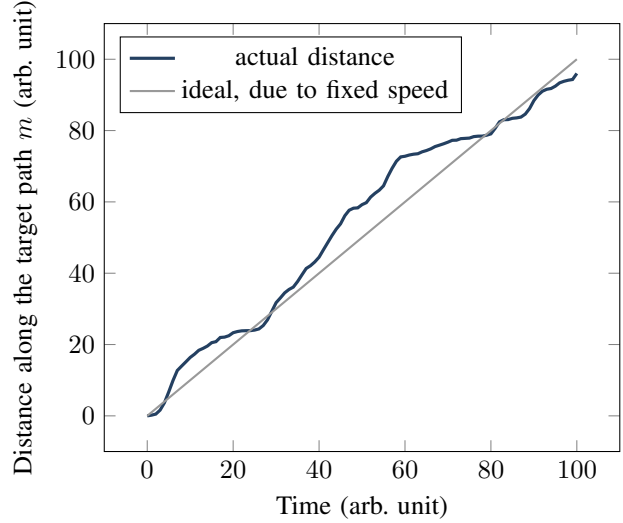
\begin{figure}[!htbp]
    \centering
    \begin{tikzpicture}
        \begin{axis}
        [
            xlabel                  = {Time (arb. unit)},
            ylabel                  = {Distance along the target path \(m\) (arb. unit)},
            cycle list name         = mycyclelist,
            legend pos              = north west,
        ]
            \addplot+ [select coords between index={0}{100}, very thick] table [ x = x, y = y ]%
            {Plotdata_linear-demo.txt};%
            \addlegendentry{actual distance};
            \addplot+ [black!40, thick] coordinates {(0,0) (100,100)};
            \addlegendentry{ideal, due to fixed speed};
         \end{axis}
\end{tikzpicture}
    \caption{If a vehicle has travelled \(m\) distance units along a path, and should continue travelling at speed \(k\) for another \(n\) units of time, then its ideal distance along the path \(d\) is a line segment with slope \(k\) and intercept \(m\) across the temporal interval \([0, n]\). The actual distance will deviate from the ideal case, but the two behaviours will typically be close.}
    \label{fig:ideal-distance-is-a-line-segment}
\end{figure}

To make this more practical, consider that typical path-following situations in marine environments require all vehicles to follow a target path at a more or less constant speed (cf.\ the lawnmower path in Fig.~\ref{fig:lawnmower}). This \emph{``nearly constant-speed''} assumption implies an ideal distance \(d\) along the path that is increasing linearly with time, even if actual speeds may deviate from this (see Fig.~\ref{fig:ideal-distance-is-a-line-segment} for a graphical representation of this intuition).

To translate the intuition from Fig.~\ref{fig:ideal-distance-is-a-line-segment} to a data-compression strategy, we note that the ideal behaviour corresponds to uniform third-degree B-splines with domain \([0, n]\) and unit interval length with spline coefficients 
\begin{equation}
    \dpl^{\textrm{ideal}}
    =
    \left[
        \dpl_{[1]}^{\textrm{ideal}}, 
        \dotsc,
        \dpl_{[n]}^{\textrm{ideal}}
    \right]
\end{equation}
equal to the vector
\begin{equation}
    m + k
    \left[
        0, \frac{1}{3}, 1, 2, 3, \dotsc, n-1, n- \frac{1}{3}, n
    \right]
    \label{equ:values-of-the-ideal-coefficients}
\end{equation}
(see the Appendix for more mathematical details behind this fact). In other words, the \(n\)-dimensional vector \(\dpl^{\textrm{ideal}}\) of coefficients corresponding to the ``ideal speed condition'' have a specific structure, and is defined by a common term \(m\) (equal to \(\dpl_[1]\)), and a specific pattern in the remaining coefficients that can be highlighted with a ``telescoping series'' operation, i.e., by computing the difference between the vector formed by the last $n-1$ coefficients and the one formed by the first $n-1$ coefficients. In other words, using the notation \(x_{[a:b]}\) to access elements \(a\) through \(b\) of \(x\),
\begin{equation}
    \label{eq:deltaspl}
    \Delta\dpl^{\textrm{ideal}}
    =
    \dpl_{[2:n]}^{\textrm{ideal}}
    -
    \dpl_{[1:n-1]}^{\textrm{ideal}}
    =
    k \left[
        \frac{1}{3}, \frac{2}{3}, \frac{3}{3}, \frac{3}{3}, \dotsc, \frac{3}{3}, \frac{2}{3}, \frac{1}{3}
    \right],
\end{equation}
that highlights that the set of differential coefficients obtained from the ``ideal constant speed'' spline is \(\{1/3,2/3,3/3\}\).

Agents need to exchange coefficients that will typically differ from these ones, since splines produced when iterating a DMPC algorithm will produce results that are different, even if close to the ideal ones. E.g., ignoring \(k\) (since the target speed is known to the agents \emph{a priori} and needs not be exchanged), a plausible vector of actual coefficients that a DMPC iteration may produce may lead to a ``telescoped vector'' equal to
\begin{equation}
    \Delta\dpl
    =
    \left[
        \frac{1.01}{3}, \frac{1.99}{3}, \frac{3.12}{3}, \frac{3.09}{3}, \dotsc, \frac{2.97}{3}, \frac{1.99}{3}, \frac{0.98}{3}
    \right] .
\end{equation}
But then, since \(\Delta\dpl^{\textrm{ideal}}\) is also prior knowledge to all agents, the information to be sent across the network is in this case
\begin{equation}
    \begin{array}{l}
    \displaystyle \frac{3}{k}
    \left( \Delta\dpl - \Delta\dpl^{\textrm{ideal}} \right)
    = \\
    \displaystyle \qquad
    \left[
        0.01, -0.01, 0.12, 0.09, \dotsc, -0.03, -0.01, -0.02
    \right] ,
    \end{array}
\label{equ:to-be-encoded-vector-of-coefficients}
\end{equation}
that shall be quantised first before being encoded in the payloads of the to-be-exchanged acoustic packets. 

\subsection{Which quantisation scheme shall be used to quantise~\eqref{equ:to-be-encoded-vector-of-coefficients}?}

Intuitively, the higher the precision used to represent such coefficients, the higher the data rate that will be required to exchange such information, but also the better the performance of the overall control scheme. Moreover, given a specific number of bits that shall be used to encode vectors as in~\eqref{equ:to-be-encoded-vector-of-coefficients}, one may also think to optimise the associated quantisation scheme, intuiting that the better the quantiser, the better the closed-loop performance.

However, this would require collecting statistics about the empirical distributions of such coefficients, and these distributions may depend on the mission (formation to maintain and path to follow, for instance). For the sake of offering results that are as general as possible, we opt for a uniform quantisation scheme. In other words, we will utilise an alphabet defined by \(\mathbb{X}_{i,f}\) for our data. Here, \(\mathbb{X}_{i,f}\) denotes the set of signed fixed-point numbers with \(i\) bits for the integer part, and \(f\) bits for the fraction part. To give an example, \(\mathbb{X}_{3,4}\) is defined by the ``strings to scalars'' look-up tables
\begin{equation}
    \left\lbrace
        \begin{array}{ccc}
            \texttt{000} & \rightarrow & 0 \\
            \texttt{001} & \rightarrow & 1 \\
            \texttt{010} & \rightarrow & 2 \\
            \texttt{011} & \rightarrow & 3 \\
            \texttt{100} & \rightarrow & -4 \\
            \texttt{101} & \rightarrow & -3 \\
            \texttt{110} & \rightarrow & -2 \\
            \texttt{110} & \rightarrow & -1 \\
        \end{array}
    \right\rbrace
    \; \cup \;
    \left\lbrace
        \begin{array}{ccc}
            \texttt{0000}, & \rightarrow & 0 \\
            \texttt{0001}, & \rightarrow & 1/16 \\
            \vdots \\
            \texttt{1111} & \rightarrow & 15/16 \\
        \end{array}
    \right\rbrace
\end{equation}
so that \texttt{010 0011} in \(\mathbb{X}_{3,4}\) encodes \(2 + 3/16 = 2.1875\).

\subsection{On the effect of choosing different quantisation levels to represent a given spline}

To help interpret the results that will be given in Section~\ref{sec:results}, we plot in Fig.~\ref{fig:effects-of-quantization-on-spline} an example of how choosing different parameters \(i\) and \(f\) to define \(\mathbb{X}_{i,f}\) translates into different levels of distortions on the quantised spline. For practical purposes, already using \(i = 3\) and \(f = 8\) provides excellent representations of the quantities that the agents shall exchange.

\begin{figure}[!htbp]
    \centering
    \tikzset{every mark/.append style={scale=0.6},}
\begin{tikzpicture}
        \begin{axis}
        [
            xlabel                  = {Time (arb. units)},
            ylabel                  = {Distance along the target path (arb. unit)},
            cycle list name         = \expandedprec cyclelist,
            mark size               = 2pt,
            legend pos              = north west
        ]
        \addplot+ [black!40, line width=2.2pt, behind path] table [ x=t, y=prec64 ] {Plotdata_fixed-point_fi-1d-spline.txt};
        \addlegendentry{double-precision}
            \foreach \myprec in {3,5,8}
                {
                \edef\addseries{\noexpand%
               	\addplot+ table [ x = t, y = prec\myprec ]%
    			{Plotdata_fixed-point_fi-1d-spline.txt};%
    			\noexpand\addlegendentry{\myprec{} fraction bits};
                }
                \addseries
    		}
	    \end{axis}
\end{tikzpicture}
    \caption{An example of how different precision levels on representing the coefficients of a given spline translate into a perturbation of such a spline. Here \(i = 3\).}
    \label{fig:effects-of-quantization-on-spline}
\end{figure}
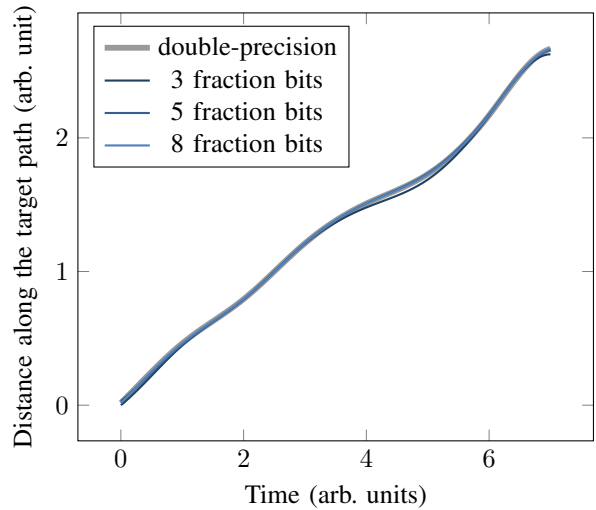

\begin{figure}[!htbp]
    \centering
    \begin{tikzpicture}
        \begin{axis}
        [
            xlabel                  = {Time (arb. units)},
            ylabel                  = {Discretization error (arb. unit)},
            xmin                    = 0,
            xmax                    = 7,
            ymin                    = -110,
            ymax                    = 10,
            ytick                   = {-100,-60,-20},
            yticklabels             = {\(10^{-5}\),\(10^{-3}\),\(10^{-1}\)},
            cycle list name         = \expandedprec cyclelist,
            mark size               = 2pt,
            legend columns          = 1,
        ]
            \foreach \myprec in {3,5,8}
                {
                \edef\addseries{\noexpand%
               	\addplot+ table [ x = t, y = prec\myprec ]%
    			{Plotdata_fixed-point_fi-1d-error-db.txt};%
    			\noexpand\addlegendentry{\myprec{} fraction bits};
                }
                \addseries
    		}
	    \end{axis}
\end{tikzpicture}
    \caption{An alternative representation of the information within Fig.~\ref{fig:effects-of-quantization-on-spline}, highlighting how the relative error quickly vanishes as one increases the number of bits used to represent the various coefficients of the splines.}
    \label{fig:error-effects-of-quantization-on-spline}
\end{figure}
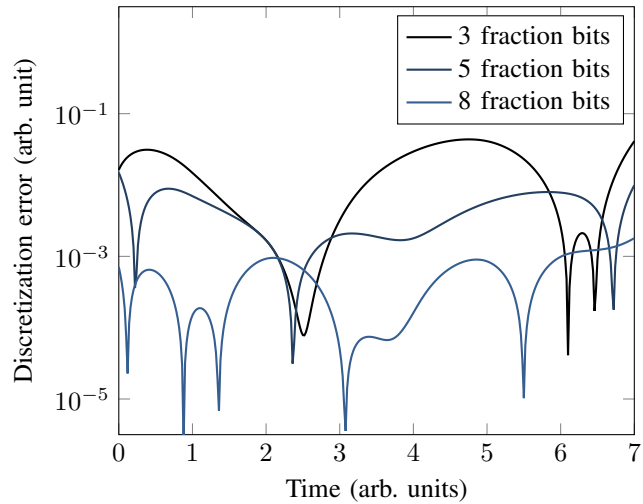

\section{Optimising the information exchange process}
\label{sec:optimising-parameters}

Independently of whether one uses a TDMA or FDMA approach, one should strive for minimising data rates subject to satisfying some given bounds on the expected control performance indexes. Assuming that the to-be-encoded coefficients (like the ones in~\eqref{equ:to-be-encoded-vector-of-coefficients}) are encoded using the alphabet \(\mathbb{X}_{i,f}\), that \(\fip\) is the number of bits defining \(\mathbb{X}_{i,f}\), that \(\fwidth\) is the number of bits used to encode the value \(m\) in~\eqref{equ:values-of-the-ideal-coefficients} (\(m\) is not known \emph{a priori}, and shall be exchanged too), and that \(\hp\) is the number of intervals defining such a spline and thus the number of coefficients in \(\Delta\dpl\), encoding a specific \(\wpl_{BA}\) as in~\eqref{eq:packets} means using
\begin{equation}
    \fwidth + \fip
    \left(
        \hp+2
    \right)
\end{equation}
bits of the available packet payload. Assume that \(\Nsize\) is the cardinality of the network, and that the generic agent \(A\) shall communicate, at each communication iteration, information about \(\Nsize-1\) splines (intuitively, one \(\wpl_{BA}\) as in~\eqref{eq:packets} for each agent \(B\) present in the network) plus the transmitter's ID. Under these assumptions the number of bits allocated in the payload will be given by
\begin{equation}
    \label{eq:psz}
    \psz
    =
    \left(
        \Nsize-1
    \right)
    \Big(
        \fwidth + \fip
        \left(
            \hp+2
        \right)
    \Big)
    +
    \lceil \log_2 \Nsize \rceil.
\end{equation}
Depending on which protocol and coding scheme one uses (possibly including the number of parity bits, size of the preamble, and so on), the \psz requirement will thus translate into a data rate \bw. That \bw should then guarantee (at least statistically) some maximum absolute error \(\msep\) and \(\mses\) on the relative positions and speeds of the vehicles, while guaranteeing enough time for the agents to transmit their packets.

It is paramount to note that the errors \(\msep\) and \(\mses\) will in the field depend on several factors, such as the actual channel conditions and the agents' positions at sea -- these factors will indeed affect both the probability of packet loss and information delay. Assuming to know an expression for the data rate \bw as a function of the parameters in the adaptation schemes proposed in Sections~\ref{sec:adapting-to-broadcast} to~\ref{sec:quantizing-splines}, one may try to optimise such adaptations by solving the program
\begin{align}
\label{eq:optproblem}
    \operatorname{minimise}\quad & \bw(\theta) \operatorname{ w.r.t.\ } \theta \in \Theta \\
    \label{eq:msepmax}
    \operatorname{subject~to}\quad & \msep(\theta) \leq \msep^\text{max} \\
    \label{eq:msesmax}
                              & \mses(\theta) \leq \mses^\text{max}
\end{align}
with \(\theta\) a placeholder for all the decision variables defined until now, and \(\Theta\) a meaningful domain for that vector. Problem~\eqref{eq:optproblem} is very hard to solve, because to the best of our knowledge, \msep and \mses may be evaluated only through simulation. They cannot be expressed analytically as functions of \(\theta\), and shall be treated as random variables (since affected by random packet loss) if one strives for field applicability. In this case~\ref{eq:msepmax} and~\ref{eq:msesmax} become inequalities in expected value, and \(\bw(\theta)\) itself cannot easily be evaluated. 

The practical need is though to have some numerical tools that may help plan missions, and decide which communication strategy to use. To this aim, we propose to relax~\eqref{eq:optproblem}, perform some numerically tractable simulation (trying to be as independent as possible from specific protocol and coding schemes), and produce repeatable results that may be used for decision-making. To this end, we thus assume a simplified model for the data rate (Section~\ref{ssec:modelling-data-rate}), and design an opportune sensitivity analysis strategy to evaluate the performance indexes \(\msep\) and \(\mses\) (Section~\ref{ssec:sensitivity-analysis}).

\subsection{Modelling the data rate associated to the proposed adaptations}
\label{ssec:modelling-data-rate}

As for modelling \bw, we assume that besides the parameters define before, any given agent broadcasts a packet every \(\Delta T\) seconds (independently of using TDMA or FDMA). This quantity is itself defined by two other parameters: the sampling time \(\Ts\) of the overall DMPC scheme (i.e., how often the MPC scheme shall run), and the number of ADMM iterations per MPC step \(\niterstep\), via the expression
\begin{equation}
    \Delta T
    =
    \frac{\Ts}{\niterstep} .
\end{equation}
Note that to allow the agents time to transmit their data, the sum of the propagation delay and the packet overhead time (e.g., its preamble) must be smaller than this \(\Delta T\). In networks with small graph diameters, the allowed time is \(\Delta T/\Nsize\) under TDMA, and \(\Delta T\) under FDMA. Given these assumptions, this means modelling the data rate for the single agent as \(\bw = \psz / \Delta T\).

\subsection{A strategy to analyse the sensitivity of the proposed adaptation mechanisms}
\label{ssec:sensitivity-analysis}

To make the proposed adaptations as effective as possible, one would like to:
\begin{itemize}
    \item on the one hand, keep the data rate requirements \bw as low as possible,
    \item on the other hand, increase the precision used when quantising the splines (hence increasing \fwidth and \fip, possibly increasing the prediction horizon \hp), and make the control loop run as often as possible (hence increasing the number of ADMM steps \niterstep and decreasing \Ts).
\end{itemize}
In short, too high compression loss, too few iterations per prediction step, too long time between prediction steps, and too few steps in the prediction horizon may increase the positional error \msep and path-tangential speed error \mses.

To understand how to tune these variables, we propose an ad-hoc sensitivity-analysis strategy to inspect how each variable influences \msep and \mses. We select a set of values for each independent variable to tune, except the width of the first coefficient \(m\), which we kept at 16~bits to provide sufficient but not excessively high precision. Because assessing \msep and \mses requires simulations, which again depend on several factors on top of the design parameters (e.g., mission types), we consider a coordinate descent approach for optimising the various parameters, outlined in Algorithm~\ref{alg:hyperparams}.

\begin{algorithm}
\caption{Ad-hoc coordinate-descent based parameters tuning}\label{alg:hyperparams}
\KwData{Collection of parameters \(\mathcal{H}\)}
\KwData{Initial values of parameters \(\mathcal{H}_0\)}
    \Begin{
        \For{\(\theta\in \mathcal{H}\)}{
            Define suitable collection \(\widehat{\theta}\) of values for \(\theta\)\;
            \For{\(\theta_i\in\widehat{\theta}\)}{
                \msep, \mses = simulate(\(\mathcal{H}_0, \theta=\theta_i\))\;
                \If{\msep and \mses acceptable}{
                    Remember \(\theta_i\)\;
                }
            }
        }
        Choose next \(\mathcal{H}_0\) among remembered values of all \(\theta\)\;
        \emph{Optional}: Save \(\mathcal{H}_0\) as parameter set\;
        \If{user believes that more data rate \bw can be saved}{
            Repeat from beginning\;
        }
    }
\end{algorithm}

In other words, starting from an initial guess of parameters values, we run a simulation that finds one series in \msep and \mses and an estimate of the required \bw. We then vary one of the parameters at the time to acquire a set of \msep and \mses series (one series per different value of the variable being explored at that step), and select a suitable next value for that variable. We then change which parameter variable is being optimised, and repeat the procedure above until no more data rate could be conserved without allowing too large \msep or \mses.

The parameters found in this way are then considered as the starting points for the simulations considered in the next section.

\section{Results}
\label{sec:results}

What we have been proposing is a set of schemes to robustify existing DMPC formation control schemes that assume ideal channels (synchronous, bidirectional, reliable) into versions that are more suitable for field deployment. Adding the robustification mechanisms will though lead to new simulated control performances, that may seem degraded with respect to the originally simulated ones (i.e., an increased deviation from target position and speed in the formation of the original DMPC). We here showcase how much these deviations may be, taking the specific case of robustifying~\cite{Matous22dMPC4FPF}; more specifically, we do this in three steps.

In Section~\ref{ssec:results-losses-wrt-ideal-comms}, we assume that no packets are lost, nor that there is need for data compression (i.e., unlimited data rate). This means checking how the adaptations alone affect performance when the communication capabilities are equivalent to that of the original DMPC algorithm.

In Section~\ref{ssec:results-sensitivity}, we investigate how sensitive the control performance is to each design parameter of the adapted DMPC algorithm and the compression strategy devised in Section~\ref{sec:quantizing-splines}. In practice, we assess for which parameters the adaptations ``break'' the controller.

In Section~\ref{ssec:results-different-parameters}, we finally compute how big a data rate is needed to maintain an arbitrarily defined acceptable tracking performance (i.e., a position error \msep that is bounded by 5\% of the formation radius, and a path-tangential speed error \mses of at most 5~cm/s in absolute value). In practice, we assess how difficult could the channel become, in terms of packet loss ratio, before tracking performance degrades too much.

\subsection{Use cases}
\label{ssec:results-use-cases}

As a first use case, we consider seabed mapping by four vehicles travelling along a lawnmower path at a constant depth, with the parameterisation of the path given by \(l=75~\mathrm{m}\) and \(r=15~\mathrm{m}\), with the vehicles ideally keeping a formation made by a 20~m wide line perpendicular to the straight sections of the lawnmower path (thus not turning with the path), with barycentre of the network on this path.

As a second use case, we consider the inspection of a vertically suspended straight pipe by six vehicles travelling in a helix path, in the sense that the formation follows a helix of radius \(20~\mathrm{m}\) (i.e., looking from ``above'' the agent moves along a circle of that radius) and spatial revolution period \(\nu_z=200\) (i.e., it takes 200 meters vertically for that agent to complete a revolution along the circle mentioned before). The six agents shall ideally form a horizontally displaced octahedron.

\subsection{Performance losses induced by the proposed adaptations w.r.t.\ the ideal communication case}
\label{ssec:results-losses-wrt-ideal-comms}

The original DMPC algorithm in~\cite{Matous22dMPC4FPF} achieves a certain formation path following performance while considering ideal communication (i.e., lossless, bidirectional, synchronous packets exchange, and double precision representations of the splines within the packets).

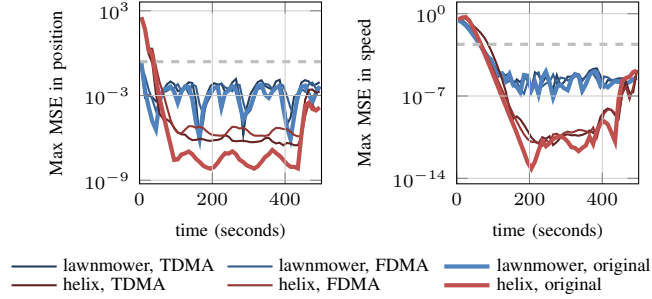
\begin{figure}[!htbp]
    \centering
    \figurefilename{tikz-compare-new-to-old}

\begin{tikzpicture}[font={\scriptsize}]
    \def\thresarray{{"0.25","0.0025"}}
    \foreach \myvariable/\myvariablelabel [count=\mycol, evaluate={\myshift=130*\mycol; \mythres=\thresarray[\mycol-1];}] in {p/position,s/speed}
    {
        \begin{axis}
        [
            ymode					= log,
    	xmode					= linear,
            xshift                  = \myshift,
            ylabel                  = {Max MSE in \myvariablelabel{}},
            legend to name          = legend \myvariable,
            legend columns          = 3,
            scale only axis,
            xmax                    = 500,
            legend style			=
		{
			draw				= none,
			fill				= none,
			inner xsep			= 0.1cm,
			inner ysep			= 0.4pt,
			nodes				= {inner ysep = 0.4pt, text width = 2.6cm}, 
			cells				= {anchor = west},
		},
        ]
            \foreach \mypath/\mypathlabel in {lawnmower/lawnmower,spiral/helix}
                {
                \edef\addseries{\noexpand%
               	\addplot table [ x = Time, y = TDMA ]%
    			{Plotdata_new-vs-old_\mypath-path-\myvariable-new-vs-old-otf.txt};%
    			\noexpand\addlegendentry{\mypathlabel, TDMA};
                \noexpand\addplot table [ x = Time, y = FDMA ]%
    			{Plotdata_new-vs-old_\mypath-path-\myvariable-new-vs-old-otf.txt};%
    			\noexpand\addlegendentry{\mypathlabel, FDMA};
    			\noexpand\addplot table [ x = Time, y = Original ]
    			{Plotdata_new-vs-old_\mypath-path-\myvariable-new-vs-old-otf.txt};
    			\noexpand\addlegendentry{\mypathlabel, original};
                }
                \addseries
    		}
            \addplot [error limit line] coordinates {(0,\mythres) (500, \mythres)};
	    \end{axis}
    }
\end{tikzpicture}

\begin{tikzpicture}[font=\scriptsize]
    \ref{legend p}
\end{tikzpicture}
    \caption{Examples of which performance losses may be introduced by the adaptations proposed in Sections~\ref{sec:adapting-to-broadcast} and~\ref{sec:adapting-to-asynch} 
    w.r.t.\ the ideal communication case.}
    \label{fig:compare-new-to-old}
\end{figure}

To quantify the losses that the adaptations proposed in Sections~\ref{sec:adapting-to-broadcast} (adapting to broadcast schemes) and~\ref{sec:adapting-to-asynch} (adapting to asynchronous communications) bring to the table, we simulate both the original DMPC algorithm and the partially adapted one for the ``reliable communication case'' (this is a requirement, since the original DMPC algorithm cannot work with unreliable communication). We then plot these performance indexes as time-series in Fig.~\ref{fig:compare-new-to-old} (note that the mean-square error values are here computed by averaging the square of respectively \msep and \mses over contiguous and non-overlapping 10-second intervals).

Fig.~\ref{fig:compare-new-to-old} plots the \msep and \mses performance indexes of the worst-performing vehicle in both use cases. We can infer that after an initial transient the additional performance losses due to the proposed adaptations are small (especially the position error in the seabed mapping application, that is practically unaffected). 

Of interest are the spikes in the position error for the seabed mapping application, which we found being an effect of approximating the sinusoidal parts of the curved parts in the lawnmower path (indeed the lawnmower curve is only once differentiable with respect to its path parameter, so a particle travelling along this path at constant speed will experience a sudden change in acceleration where the straights meet the turns).

\subsection{Sensitivity analysis with respect to the defining parameters}
\label{ssec:results-sensitivity}

\begin{figure*}[htbp]
    \centering
    \begin{tikzpicture}[font=\scriptsize]
    \node [align=center, font=\normalsize, anchor=base] at (0.7*\twocolumnwidth, 0.6) {Position error};
    \node [align=center, font=\normalsize, anchor=base] at (1.7*\twocolumnwidth, 0.6) {Speed error};
    \foreach \x in {0,1} {
        \node [align=center, font=\small, anchor=base] at (\x.525*\twocolumnwidth, 0.2) {Helix};
        \node [align=center, font=\small, anchor=base] at (\x.875*\twocolumnwidth, 0.2) {Lawnmower};
    }
    \def\hparray{3,4,6,10}
    \def\tsarray{1,3,6,9}
    \def\qsarray{10,47,220,1000,2200,4700,10000}
    \def\niterarray{1,2,3,4}
    \def\precarray{2,4,8,12}
    \def\decimarray{1,2,3,4}
    \def\thresarray{{"0.25","0.0025"}}
    \def\yticktenposition{"(100,1e-8) (100,1e-4) (100, 1)"}
    \def\yticktenspeed{"(100,1e-12) (100,1e-6) (100, 1)"}
    \def\whichvariable{{"$h_p$","$T_s$","$N_\text{ips}$","$w_\text{fi}$","$N_d$"}}
    \def\witherrorline{1}
    \foreach \mypath [count=\mycol] in {helix,lawnmower}
    \foreach \mysensitivity [count=\myrow, evaluate={\myyshift=-0.35*\myrow*\twocolumnwidth;}] in {hp,ts,niter,prec}
    {
        \foreach \myvariablelabel [count=\mycolumn, evaluate={\myshift=0.35*\twocolumnwidth*\mycol+1*\twocolumnwidth*(\mycolumn-1);\mytens=ifthenelse(\mycolumn==1,\yticktenposition,\yticktenspeed);\mythres=\thresarray[\mycolumn-1];}] in {position,speed}
        {
            \begin{axis}
            	[
                    width=0.35\twocolumnwidth,
                    height=0.35\twocolumnwidth,
                    ymode					= log,
            		xmode					= linear,
                    xmax                    = 500,
                    xtick style             = {draw opacity={ifthenelse(\myrow==4,1,0)}},
                    ymin                    = {ifthenelse(\mycolumn==1,1e-10,1e-13)},
                    ymax                    = {ifthenelse(\mycolumn==1,1e2,10)},
                    ytick                   = data,
                    xtick                   = {0, 200, 400},
                    xticklabels             =  {\ifthenelse{\equal{\myrow}{4}}{0}{\empty},\ifthenelse{\equal{\myrow}{4}}{200}{\empty},\ifthenelse{\equal{\myrow}{4}}{400}{\empty}},
                    ytick style             = {draw opacity={ifthenelse(\mycol==1,1,0)}},
                    yticklabels              = {\ifthenelse{\equal{\mycol}{1}}{\ifthenelse{\equal{\mycolumn}{1}}{\(10^{-8}\)}{\(10^{-12}\)}}{\empty},\ifthenelse{\equal{\mycol}{1}}{\ifthenelse{\equal{\mycolumn}{1}}{\(10^{-4}\)}{\(10^{-6}\)}}{\empty},\ifthenelse{\equal{\mycol}{1}}{\(10^0\)}{\empty}},
                    xlabel                  =\ifthenelse{\equal{\myrow}{4}}{time (seconds)}{},
                    xshift                  = \myshift,
                    yshift                  = \myyshift,
                    ylabel                  = \ifthenelse{\equal{\mycol}{1}}{max MSE in \myvariablelabel}{},
                    ylabel style={align=center, text width=2cm},
                    legend columns          = 1,
                    scale only axis,
                    cycle list name         =\expanded\mysensitivity cyclelist
                ]
                \edef\cheattheyticks{\noexpand\addplot+ [draw opacity=0] %
                    coordinates {\mytens};
                }
                \cheattheyticks
                \addlegendentry{};
                \edef\myaddtitle{\noexpand%
                \ifthenelse{\noexpand\equal{\mypath\myvariablelabel}{lawnmowerspeed}}{
                \noexpand\addlegendimage{empty legend};
                \noexpand\pgfmathsetmacro{\noexpand\myvar}{\noexpand\whichvariable[\myrow-1]}
                \noexpand\addlegendentry{\noexpand\myvar};
                }{}
                }
                \myaddtitle
                \edef\thisarray{\csname\mysensitivity array\endcsname}
                \foreach \ts in \thisarray {
                    \edef\addseries{\noexpand%
                    \addplot table [ x = Time, y = \mysensitivity\ts ]%
                    {Plotdata_tdma-\mypath-sensitivity_sensitivity-\mypath-tdma-\myvariablelabel.txt};%
                    \noexpand\ifthenelse{\noexpand\equal{\mypath\myvariablelabel}{lawnmowerspeed}}{%
                    \noexpand\addlegendentry{$\ts$};}{}
                    }
                    \addseries
                }
                \addplot+ [error limit line, dashed, very thick, draw opacity=\witherrorline] coordinates {(0,\mythres) (500, \mythres)};
            \end{axis}
        }
    };
\end{tikzpicture}
    \caption{The array of plots shows the performance indexes mean-square position error \Msep and mean-square speed error \Mses when the DMPC algorithm is run using TDMA. The default parameters are used for all variables but one, which changed between simulations. The default values are \(\hp=10, \Ts=5, \niterstep=5, \fip=53\). The same variable was changed in all plots that are on the same row. Which variable was changed, and to what, is shown in the legend to the far right. The plots in the left half of the figure shows series of \Msep, and the plots in the right half show series of \Mses. Within each ``performance-index half'', the left column shows results from the use case with the helix path, and the right column shows results from the use case with the lawnmower path.}
    \label{fig:allerror-start}
\end{figure*}
\begin{figure*}[htbp]
    \centering
    \begin{tikzpicture}[font=\scriptsize]
    \node [align=center, font=\normalsize, anchor=base] at (0.7*\twocolumnwidth, 0.6) {Position error};
    \node [align=center, font=\normalsize, anchor=base] at (1.7*\twocolumnwidth, 0.6) {Speed error};
    \foreach \x in {0,1} {
        \node [align=center, font=\small, anchor=base] at (\x.525*\twocolumnwidth, 0.2) {Helix};
        \node [align=center, font=\small, anchor=base] at (\x.875*\twocolumnwidth, 0.2) {Lawnmower};
    }
    \def\hparray{3,4,6,10}
    \def\tsarray{1,3,6,9}
    \def\qsarray{10,47,220,1000,2200,4700,10000}
    \def\niterarray{1,2,3,4}
    \def\precarray{2,4,8,12}
    \def\decimarray{1,2,3,4}
    \def\thresarray{{"0.25","0.0025"}}
    \def\yticktenposition{"(100,1e-8) (100,1e-4) (100, 1)"}
    \def\yticktenspeed{"(100,1e-12) (100,1e-6) (100, 1)"}
    \def\whichvariable{{"$h_p$","$T_s$","$N_\text{ips}$","$w_\text{fi}$","$N_d$"}}
    \def\witherrorline{1}
    \foreach \mypath [count=\mycol] in {helix,lawnmower}
    \foreach \mysensitivity [count=\myrow, evaluate={\myyshift=-0.35*\myrow*\twocolumnwidth;}] in {hp,ts,niter,prec}
    {
        \foreach \myvariablelabel [count=\mycolumn, evaluate={\myshift=0.35*\twocolumnwidth*\mycol+1*\twocolumnwidth*(\mycolumn-1);\mytens=ifthenelse(\mycolumn==1,\yticktenposition,\yticktenspeed);\mythres=\thresarray[\mycolumn-1];}] in {position,speed}
        {
            \begin{axis}
            	[
                    width=0.35\twocolumnwidth,
                    height=0.35\twocolumnwidth,
                    ymode					= log,
            		xmode					= linear,
                    xmax                    = 500,
                    xtick style             = {draw opacity={ifthenelse(\myrow==4,1,0)}},
                    ymin                    = {ifthenelse(\mycolumn==1,1e-10,1e-13)},
                    ymax                    = {ifthenelse(\mycolumn==1,1e2,10)},
                    ytick                   = data,
                    xtick                   = {0, 200, 400},
                    xticklabels             =  {\ifthenelse{\equal{\myrow}{4}}{0}{\empty},\ifthenelse{\equal{\myrow}{4}}{200}{\empty},\ifthenelse{\equal{\myrow}{4}}{400}{\empty}},
                    ytick style             = {draw opacity={ifthenelse(\mycol==1,1,0)}},
                    yticklabels              = {\ifthenelse{\equal{\mycol}{1}}{\ifthenelse{\equal{\mycolumn}{1}}{\(10^{-8}\)}{\(10^{-12}\)}}{\empty},\ifthenelse{\equal{\mycol}{1}}{\ifthenelse{\equal{\mycolumn}{1}}{\(10^{-4}\)}{\(10^{-6}\)}}{\empty},\ifthenelse{\equal{\mycol}{1}}{\(10^0\)}{\empty}},
                    xshift                  = \myshift,
                    yshift                  = \myyshift,
                    xlabel                  =\ifthenelse{\equal{\myrow}{4}}{time (seconds)}{},
                    ylabel                  = \ifthenelse{\equal{\mycol}{1}}{max MSE in \myvariablelabel}{},
                    ylabel style={align=center, text width=2cm},
                    legend columns          = 1,
                    scale only axis,
                    cycle list name         =\expanded\mysensitivity cyclelist
                ]
                \edef\cheattheyticks{\noexpand\addplot+ [draw opacity=0] %
                    coordinates {\mytens};
                }
                \cheattheyticks
                \addlegendentry{};
                \edef\myaddtitle{\noexpand%
                \ifthenelse{\noexpand\equal{\mypath\myvariablelabel}{lawnmowerspeed}}{
                \noexpand\addlegendimage{empty legend};
                \noexpand\pgfmathsetmacro{\noexpand\myvar}{\noexpand\whichvariable[\myrow-1]}
                \noexpand\addlegendentry{\noexpand\myvar};
                }{}
                }
                \myaddtitle
                \edef\thisarray{\csname\mysensitivity array\endcsname}
                \foreach \ts in \thisarray {
                    \edef\addseries{\noexpand%
                    \addplot table [ x = Time, y = \mysensitivity\ts ]%
                    {Plotdata_fdma-\mypath-sensitivity_sensitivity-\mypath-fdma-\myvariablelabel.txt};%
                    \noexpand\ifthenelse{\noexpand\equal{\mypath\myvariablelabel}{lawnmowerspeed}}{%
                    \noexpand\addlegendentry{$\ts$};}{}
                    }
                    \addseries
                }
                \addplot+ [error limit line, draw opacity=\witherrorline] coordinates {(0,\mythres) (500, \mythres)};
            \end{axis}
        }
    };
\end{tikzpicture}
    \caption{The array of plots shows the performance indexes mean-square position error \Msep and mean-square speed error \Mses when the DMPC algorithm is run using FDMA. The default parameters are used for all variables but one, which changed between simulations. The default values are \(\hp=10, \Ts=5, \niterstep=5, \fip=53\). The same variable was changed in all plots that are on the same row. Which variable was changed, and to what, is shown in the legend to the far right. The plots in the left half of the figure shows series of \Msep, and the plots in the right half show series of \Mses. Within each ``performance-index half'', the left column shows results from the use case with the helix path, and the right column shows results from the use case with the lawnmower path.}
    \label{fig:allerror-fdma}
\end{figure*}

Compressing data inevitably leads to precision loss, as discussed in Section~\ref{sec:quantizing-splines}. MPC parameters (horizon length \hp and sampling time \Ts) will also affect such performance. We here quantitatively assess these dependencies.

Fig.s~\ref{fig:allerror-start} and~\ref{fig:allerror-fdma} show a collection of simulation results for reliable communication (no instances of lost packets), in the sense that each subfigure shows how the mean-square errors \Mses in position and \Msep in speed vary with each single parameter (Fig.~\ref{fig:allerror-start} referring then to TDMA schemes, and~\ref{fig:allerror-fdma} to FDMA ones). More specifically, each row of these tables correspond to a specific parameter, while each column refers to a specific control performance index (\Mses or \Msep) and specific use case path (lawnmower or helix). In all these simulations the precision used for representing the first coefficient of the splines \fwidth was kept at 32~bits, so to let the agents represent a position up to 1000~km along the path with at least decimetre precision, which we believe is sufficient for any practical marine application. Analysing the results means then inspecting per-row how the different parameters affect the performance indexes.

From the first rows of the two figures we note that increasing the MPC horizon length \hp raises the data rate requirements almost proportionally (since this implies having longer representations for the splines), but also improves the overall control performance (as expected in general MPC situations). There seems to be the need then for a sufficiently long receding-horizon length \hp to attain a practically stable controller both in the TDMA and FDMA cases, and this means that the requirements for how many splines coefficients shall be sent by the agents has a lower limit.

From the second rows of the two figures we note that the MPC sampling period \Ts has a significant impact on performance, both for TDMA and for FDMA. Intuitively, increasing \Ts lowers the data rate requirements owing to the reciprocal relationship in the definition of \(\bw=\psz\niterstep/\Ts\), but at the same time, increasing \Ts corresponds to lowering the frequency of the control loop iterations. This clearly raises both positional and speed errors, and slows down how fast the controller may reject disturbances. Results seem to indicate that, at least for the scenarios considered, \Ts shall not be much higher than 6.

From the third rows of the two figures we note that the number of ADMM iterations that are run per MPC step, i.e., \niterstep, has actually very little effect on the performance indexes \msep and \mses both for TDMA and FDMA cases (with the TDMA case suggesting that it is enough that at least two iterations are made per step to achieve good performance). This shows high savings potential, because increasing \niterstep raises the data rate requirements of the scheme proportionally.

From the fourth rows of the two figures we note that the amount of bits used to represent the various splines, i.e., \fip, has as expected significant effects on the control performance indexes for both TDMA and FDMA cases. More precisely, we note the same ``diminishing returns'' effect that we noted in Fig.s~\ref{fig:effects-of-quantization-on-spline} and~\ref{fig:error-effects-of-quantization-on-spline}: the relative performance degradation quickly vanishes as one increases the number of bits used to represent the various coefficients of the splines. One may thus safely assume that using four bits for the fractional part for each spline coefficient is enough.

\subsection{Performance with different parameter sets}
\label{ssec:results-different-parameters}

From the sensitivity analysis, we selected parameter sets that indicatively gave acceptable performance and had much lower data rate requirements of the algorithm than what a parameter set corresponding to the original DMPC algorithm would demand. The simulations used a constant and global probability of packet loss per link to enable demonstrating the efficacy of the packet-loss handler.

Table~\ref{tab:psets} shows different sets of parameters for which we calculated data rate requirements (again for a six-vehicle network using TDMA). Parameter set~0 was chosen to match the original DMPC algorithm as closely as possible. For FDMA, the per-agent data rate requirements decrease by a factor equal to the network size, while the network-wide data rate remains the same, as does the time-bandwidth allowance per agent.

\begin{table}[htbp]
    \centering
    \caption{Parameter value sets.}
    \label{tab:psets}
    \begin{tabular}{cccc}
        \toprule
        parameter & value & value & value \\
        \& symbol & in set 0 & in set 1 & in set 2 \\
        \midrule
        MPC horizon length \hp & 10 & 6 & 6 \\
        sampling time \Ts & 5 & 8 & 12 \\
        ADMM iterations per MPC step \niterstep & 5 & 3 & 2 \\
        quantisation index \fip (f. part) & 53 & 10 & 5 \\
        data rate requirements \emph{\bw} & \emph{22~458} & \emph{1627} & \emph{523} \\
        \bottomrule
    \end{tabular}
\end{table}

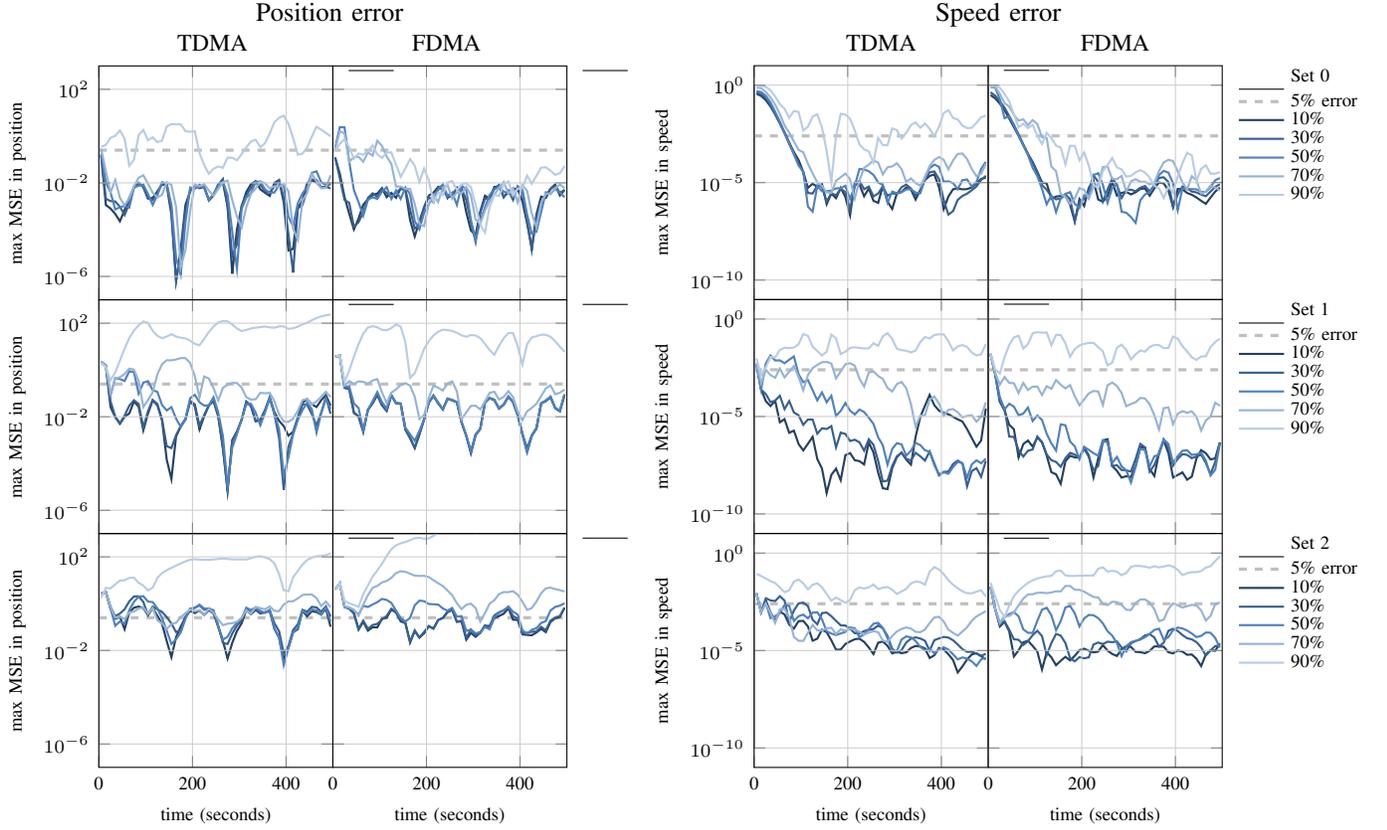
\begin{figure*}[htbp]
    \centering
    \figurefilename{tikz-sensitivity}

\begin{tikzpicture}[font=\scriptsize]
    \node [align=center, font=\normalsize, anchor=base] at (1.675*\twocolumnwidth, 0.6) {Position error};
    \node [align=center, font=\normalsize, anchor=base] at (2.675*\twocolumnwidth, 0.6) {Speed error};
    \foreach \x in {1,2} {
        \node [align=center, font=\small, anchor=base] at (\x.5*\twocolumnwidth, 0.2) {TDMA};
        \node [align=center, font=\small, anchor=base] at (\x.85*\twocolumnwidth, 0.2) {FDMA};
    }
    \def\lossarray{10,30,50,70,90} 
    \def\psetarray{0,1,2}
    \def\ylimlow{{"1e-7","1e-11"}}
    \def\ylimhigh{{1e3,10}}
    \def\ythres{{0.25,0.0025}}
    \def\yticktenposition{"(100,1e-6) (100,1e-2) (100, 100)"}
    \def\yticktenspeed{"(100,1e-10) (100,1e-5) (100, 1)"}
    \foreach \myvariable/\myvariablelabel [count=\mycol, evaluate={\myshift=\twocolumnwidth*\mycol; \myymin=\ylimlow[\mycol-1];\myymax=\ylimhigh[\mycol-1];\mythres=\ythres[\mycol-1];}] in {p/position,s/speed}
    \foreach [count=\myrow] \ts in \psetarray \foreach [count=\myprotocol,evaluate={\mytens=ifthenelse(\mycol==1,\yticktenposition,\yticktenspeed);}] \strat in {\empty,"-fdma"}
    {
    \begin{axis}
        [
            width=0.35\twocolumnwidth,
            height=0.35\twocolumnwidth,
            xtick style             = {draw opacity={ifthenelse(\myrow==3,1,0)}},
            ytick                   = data,
            xtick                   = {0, 200, 400},
            xticklabels             =  {\ifthenelse{\equal{\myrow}{3}}{0}{\empty},\ifthenelse{\equal{\myrow}{3}}{200}{\empty},\ifthenelse{\equal{\myrow}{3}}{400}{\empty}},
            xtick style             = {draw opacity={ifthenelse(\myrow==3,1,0)}},
            ytick style             = {draw opacity={ifthenelse(\myprotocol==1,1,0)}},
            yticklabels              = {\ifthenelse{\equal{\myprotocol}{1}}{\ifthenelse{\equal{\mycol}{1}}{\(10^{-6}\)}{\(10^{-10}\)}}{\empty},\ifthenelse{\equal{\myprotocol}{1}}{\ifthenelse{\equal{\mycol}{1}}{\(10^{-2}\)}{\(10^{-5}\)}}{\empty},\ifthenelse{\equal{\myprotocol}{1}}{\ifthenelse{\equal{\mycol}{1}}{\(10^{2}\)}{\(10^{0}\)}}{\empty}},
            ymode					= log,
            xmode					= linear,
            ymin                    = \myymin,
            ymax                    = \myymax,
            xmax                    = 500,
            ylabel                  = {Max MSE\\in \myvariablelabel},
            xshift                  = 0.98*\myshift + 0.35*\twocolumnwidth*\myprotocol,
            yshift                  = -0.35*\twocolumnwidth*\myrow,
            xlabel                  =\ifthenelse{\equal{\myrow}{3}}{time (seconds)}{},
            ylabel                  = \ifthenelse{\equal{\myprotocol}{1}}{max MSE in \myvariablelabel}{},
            ylabel style={align=center, text width=2cm},
            legend columns          = 1,
            scale only axis,
            cycle list name         =pset cyclelist%
        ]

        \edef\myaddtitle{\noexpand%
        \ifthenelse{\noexpand\equal{\myprotocol\myvariablelabel}{2speed}}{
        \noexpand\addlegendimage{empty legend};
        \noexpand\addlegendentry{Set~\ts};
        }{}
        }
        \myaddtitle
        \edef\cheattheyticks{\noexpand\addplot+ [draw opacity=0] %
            coordinates {\mytens};
        }
        \cheattheyticks
        \addlegendentry{};
        \addplot [error limit line] coordinates { (0, \mythres) (500, \mythres) };
        \ifthenelse{\equal{\myprotocol\myvariablelabel}{2speed}}{\addlegendentry{5\% error};}{}
        \foreach \loss in \lossarray {
            \edef\addseries{\noexpand%
            \addplot table [ x = Time, y = Set\ts Loss\loss ]%
            {Plotdata_param-sets_lawnmower-mse\myvariable\strat-allpsets.txt};%
            \noexpand\ifthenelse{\noexpand\equal{\myprotocol\myvariablelabel}{2speed}}{\noexpand\addlegendentry{\loss\noexpand\%};}{}}
            \addseries
        }
    \end{axis}
    }
\end{tikzpicture}
    \caption{The array of plots shows the performance indexes mean-square position error \Msep and mean-square speed error \Mses when the DMPC algorithm is run using the parameter sets in Table~\ref{tab:psets}, first under TDMA, then under FDMA. The packet loss ratio and the parameter set used is shown in the legend to the far right; the legends apply row-wise. The plots in the left half of the figure shows series of \Msep, and the plots in the right half show series of \Mses. Within each ``performance-index half'', the left column shows results from using TDMA, and the right column shows results from using FDMA.}
    \label{fig:mse-in-s-lawnmower}
\end{figure*}

Fig.~\ref{fig:mse-in-s-lawnmower} shows the positional and path-speed errors under the selected parameter sets, using TDMA or FDMA. The use case was simulated five times for each parameter set, setting the packet loss ratio to \(10\%, 30\%, 50\%, 70\%, 90\%\), in order. The black line marks the error threshold that we wish the network to stay below.

We see here that parameter set~1 leads to acceptable performance under TDMA, that the maximum \Msep eventually stays below 5\% (marked with a black line) of the formation radius, and that the maximum \Mses eventually falls below 5\%, even under 70\% packet loss. It also requires almost 14 times lower data rate than the original DMPC algorithm would have required, showing that our adaptations enable DMPC in an AUV network at a much lower data rate than would have been required without data compression.

\begin{figure}[htbp]
    \centering
    \begin{tikzpicture}[font=\scriptsize]
    \def\tsarray{10,11,12}
    \def\precarray{1,2,5}
    \def\thresarray{{"0.25","0.0025"}}
    \def\yticktenposition{"(100,1e-2) (100,1e-1) (100, 1)"}
    \def\yticktenspeed{"(100,1e-8) (100,1e-4) (100, 1)"}
    \def\whichvariable{{"$T_s$","$w_\text{fi}$"}}
    \def\witherrorline{1}
    \def\mycol{1}
    \def\mypath{lawnmower}
    \foreach \mysensitivity [count=\myrow, evaluate={\myyshift=-0.3*\myrow*\twocolumnwidth;}] in {ts,prec}
    {
        \foreach \myvariablelabel/\myvariableshort [count=\mycolumn, evaluate={\myshift=0.3*\twocolumnwidth*\mycol+0.55*\twocolumnwidth*(\mycolumn-1);\mytens=ifthenelse(\mycolumn==1,\yticktenposition,\yticktenspeed);\mythres=\thresarray[\mycolumn-1];}] in {position/p,speed/s}
        {
            \begin{axis}
            	[
                    width=0.3\twocolumnwidth,
                    height=0.3\twocolumnwidth,
                    ymode					= log,
            		xmode					= linear,
                    xmax                    = 500,
                    xtick style             = {draw opacity={ifthenelse(\myrow==2,1,0)}},
                    ymin                    = {ifthenelse(\mycolumn==1,1e-3,1e-9)},
                    ymax                    = {ifthenelse(\mycolumn==1,10,10)},
                    ytick                   = data,
                    xtick                   = {0, 200, 400},
                    xticklabels             =  {\ifthenelse{\equal{\myrow}{2}}{0}{\empty},\ifthenelse{\equal{\myrow}{2}}{200}{\empty},\ifthenelse{\equal{\myrow}{2}}{400}{\empty}},
                    ytick style             = {draw opacity={ifthenelse(\mycol==1,1,0)}},
                    yticklabels              = {\ifthenelse{\equal{\mycol}{1}}{\ifthenelse{\equal{\mycolumn}{1}}{\(10^{-2}\)}{\(10^{-8}\)}}{\empty},\ifthenelse{\equal{\mycol}{1}}{\ifthenelse{\equal{\mycolumn}{1}}{\(10^{-1}\)}{\(10^{-4}\)}}{\empty},\ifthenelse{\equal{\mycol}{1}}{\(10^0\)}{\empty}},
                    xshift                  = \myshift,
                    yshift                  = \myyshift,
                    xlabel                  =\ifthenelse{\equal{\myrow}{2}}{time (seconds)}{},
                    ylabel                  = \ifthenelse{\equal{\mycol}{1}}{max MSE in \myvariablelabel}{},
                    ylabel style={align=center, text width=2cm},
                    legend columns          = 1,
                    scale only axis,
                    label shift             = -5pt,
                    cycle list name         =\expanded\mysensitivity cyclelist
                ]
                \edef\cheattheyticks{\noexpand\addplot+ [draw opacity=0] %
                    coordinates {\mytens};
                }
                \cheattheyticks
                \addlegendentry{};
                \edef\myaddtitle{\noexpand%
                \ifthenelse{\noexpand\equal{\mypath\myvariablelabel}{lawnmowerspeed}}{
                \noexpand\addlegendimage{empty legend};
                \noexpand\pgfmathsetmacro{\noexpand\myvar}{\noexpand\whichvariable[\myrow-1]}
                \noexpand\addlegendentry{\noexpand\myvar};
                }{}
                }
                \myaddtitle
                \edef\thisarray{\csname\mysensitivity array\endcsname}
                \foreach \ts in \thisarray {
                    \edef\addseries{\noexpand%
                    \addplot table [ x = Time, y = \mysensitivity\ts.0 ]%
                    {Plotdata_small-tuning_show-undersample-\myvariableshort-\mysensitivity.txt};%
                    \noexpand\ifthenelse{\noexpand\equal{\mypath\myvariablelabel}{lawnmowerspeed}}{%
                    \noexpand\addlegendentry{$\ts$};}{}
                    }
                    \addseries
                }
                \addplot+ [error limit line, draw opacity=\witherrorline] coordinates {(0,\mythres) (500, \mythres)};
            \end{axis}
        }
    };
\end{tikzpicture}
    \caption{Fine-tuning parameter set~2 to meet the position-error constraint under FDMA. Left: position error. Right: speed error.}
    \label{fig:better-than-set-2}
\end{figure}
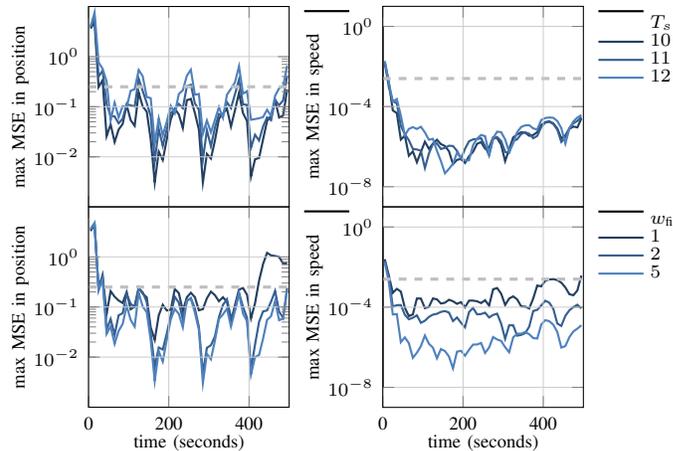

We also see that parameter set~2 meets the target \mses, but fails to meet the \msep target. A possible reason to this is that the agents do not iterate often enough in absolute time, so that the trajectory undersamples the path. The time between iterations with parameter set~2 becomes 6~s, which is more than twice as long as with parameter set~1 (2.67~s). The undersampling leads to a loss in precision; as Fig.~\ref{fig:better-than-set-2} (top row) shows, lowering \Ts to 10~s in parameter set~2, for an allowed time of 5~s between iterations, lets the vehicles attain the \msep target under FDMA. The resulting increase in data rate can be compensated by lowering the precision down to two fractional bits and still attain the \msep target, as Fig.~\ref{fig:better-than-set-2} (bottom row) shows. In fact, applying these two changes lowers \bw to 502~bit/s in a six-node network.

\section{Conclusion}
\label{sec:conclusion}

Communicating digital information acoustically means both risking not being able to deliver a packet, and needing meaningful strategies on how to quantise that information. In this paper, we explored how formation control algorithms that treat information exchange in idealistic manners can be transformed into schemes that are aware of the non-idealities of reality.

Adapting an ``ideal'' algorithm to make it practically implementable means indeed adding layers that may degrade the performance and properties that may be expected by just simulating the original control scheme. The questions revolve then around ``how much degradation'' the adaptations bring to the table, and when the constraints on the acoustic channel (essentially its bandwidth) imply unacceptable degradation levels (the latter typically being a use-case-specific quantity).

The paper thus proposes a number of mechanisms that may relax the assumptions of synchronous, bidirectional, reliable, and unlimited data rate communication. These are characterised by a set of opportune parameters (e.g., number of bits used to quantise the splines coefficients), and choosing a given set of parameters implies some expected closed loop performance indexes and some expected data rate requirements.

Interestingly, most of the data rate savings that can be achieved through the proposed modifications come from reducing the precision of the numbers encoded in the to-be-sent packets. Additional data rate savings can be obtained by lowering the frequency of the control action, and by shortening the prediction horizon of a model predictive control approach -- all these factors though up to a certain point, beyond which performance degradation spikes.

Summarising, we implicitly developed a framework for which, given a certain use-case and a constraint on the available bandwidth, one may thus search which parameters give the best closed loop performance indexes of the practically implementable version of the controller, and verify whether the performance degradation levels are acceptable or not for that specific use-case.

The framework opens also up an interesting future extensions: assuming that the parameters are capable of attaining satisfactory performance, one may think about tuning the acoustic transmission power dynamically, so that the probability of being detected by adversaries or disturbing marine life on the site decreases, while at the same time ensuring a minimum closed loop performance.

\section*{Acknowledgements}

This work has received funding from the Research Council of Norway through grant numbers 302435 and 344847.
The authors also wish to thank Josef Matou{\v{s}}, Hefeng Dong and John Potter for providing insights and suggestions about the various techniques.

\appendix
\section{Derivation of~\eqref{equ:values-of-the-ideal-coefficients}}
In this section, we derive the coefficients \spl of a uniform third-degree B-spline \Spl that equals the line segment \(y=kx+m\) for \(0\leq x \leq n\). \Spl is defined by the break points \(\{0,1,\dotsc,n\}\). The interested reader may refer to e.g.~\cite{deBoor78splines,Bojanov93splines} for further reading.

Let \(\vec{t}=(0,0,0,0,1,\dotsc,n,n,n,n)\) be the \emph{knot vector} that defines the uniform third-degree \Spl, and let \spl be a column vector. We will use the subscripted \(\Spl_{j,d}\) to indicate the \(j\)-th basis function of \Spl as a \(d\)-th degree B-spline. The basis functions are defined degree-recursively by the relationship~\cite{Lyche07bspline}
\begin{equation}
\label{eq:bspline-recurse}
    \Spl_{j,d}(x)
    =
    \frac{
        x-t_{[j]}
    }{
        t_{[j+d]}-t_{[j]}
    }
    \Spl_{j,d-1}(x)
    +
    \frac{
        t_{[j+d+1]}-x
    }{
        t_{[j+d+1]}-t_{[j+1]}
    }
    \Spl_{j,d-1}(x),
\end{equation}
where we define the basic case \(\Spl_{j,0}=u(x-t_{[j]})u(t_{[j+1]}-x)\) and \(u(x)\) as the unit step function. Graphically, \(\Spl_{j,0}\) is a rectangle equal to one over at most one interval. \(\Spl_{j,1}\) is a triangle over two intervals, or half a triangle if one of the intervals are zero width, such that the value of \(\Spl_{j,1}\) is maximum at the zero-width interval.

We omit the recursive calculations necessary from~\eqref{eq:bspline-recurse} in interest of space and present only the four possible shapes of the basis functions. Note that we assume \(n\geq 4\) to show all possible basis functions of a uniform \Spl. The recursive calculations give that
\begin{equation}
    \label{eq:bsplinebasisfunc}
    \begin{pmatrix}
        \Spl_{1,3} \\ 
        \Spl_{2,3} \\ 
        \Spl_{3,3} \\ 
        \Spl_{4,3} 
    \end{pmatrix}
    =
    \left(
    \sum_{p=0}^3A_px^p
    \right)
    \begin{pmatrix}
        u(x)u(1-x)\\
        u(x-1)u(2-x)\\
        u(x-2)u(3-x)\\
        u(x-3)u(4-x)
    \end{pmatrix},
\end{equation}
where the polynomial coefficient matrices
\begin{gather*}
A_0=
\begin{pmatrix}
1 & 0 & 0 & 0 \\
0 & 2 & 0 & 0 \\
0 & -\frac{3}{2} & \frac{9}{2} & 0 \\
0 & \frac{2}{3} & \frac{22}{3} & \frac{32}{3}
\end{pmatrix},
A_1=
\begin{pmatrix}
-3 & 0 & 0 & 0 \\
 3 & -3 & 0 & 0 \\
 0 & \frac{9}{2} & -\frac{9}{2} & 0 \\
 0 & -2 & 10 & -8
\end{pmatrix},
\\
A_2=
\begin{pmatrix}
 3 & 0 & 0 & 0 \\
 \frac{9}{2} & \frac{3}{2} & 0 & 0 \\
 \frac{3}{2} & -3 & \frac{3}{2} & 0 \\
 0 & 2 & -4 & 2
\end{pmatrix},
A_3=
\begin{pmatrix}
 -1 & 0 & 0 & 0 \\
 \frac{7}{4} & \frac{-1}{4} & 0 & 0 \\
 \frac{-11}{12} & \frac{7}{12} & \frac{-2}{12} & 0 \\
 \frac{1}{6} & \frac{-1}{2} & \frac{1}{2} & \frac{-1}{6}
\end{pmatrix}.
\end{gather*}
\(\Spl_{5,3}\) up to \(\Spl_{n,3}\) are shifted copies of \(\Spl_{4,3}\), and the last three \(\Spl_{j,3}\) are opportunely shifted and mirrored copies of the first three \(\Spl_{j,3}\). Thanks to these properties, we will only need to give the coefficients required to uniquely define the first four intervals.

First, consider the interval \([0, 1]\). From~\eqref{eq:bsplinebasisfunc}, we have that
\begin{equation}
\label{eq:bsplfirstintvl}
    u(x)u(1-x)
    \begin{pmatrix}
        \Spl_{1,3} & \Spl_{2,3} & \Spl_{3,3} & \Spl_{4,3}
    \end{pmatrix}
    \spl_{[1:4]}=u(x)u(1-x)(kx+m).
\end{equation}
Expressed in the polynomial coefficients, we have that
\begin{equation}
    \label{eq:polyfirstintvl}
    \begin{matrix}
        1\colon \\
        x\colon \\
        x^2\colon \\
        x^3\colon 
    \end{matrix}
    \begin{pmatrix}
        1 & 0 & 0 & 0 \\
        -3 & 3 & 0 & 0 \\
        3 & \frac{-9}{2} & \frac{3}{2} & 0 \\
        -1 & \frac{7}{4} & \frac{-11}{12} & \frac{1}{6} 
    \end{pmatrix}
    \spl_{[1:4]}
    =
    \begin{pmatrix}
        m \\
        k \\
        0 \\
        0
    \end{pmatrix},
\end{equation}
which has the unique solution \(\spl_{[1:4]}=m + k\left[0, \frac{1}{3}, 1, 2\right]\). \(\spl_{[p]}\) can then be calculated recursively for \(p=5,\dotsc,n+3\) by solving~\eqref{eq:polyfirstintvl}, using \(\Spl_{p,3}\) and the three preceding basis functions together with \(\spl_{[p-3:p-1]}\), eventually solving for the entire \spl as in~\eqref{equ:values-of-the-ideal-coefficients}.

\printbibliography

\end{document}